
%

\baselineskip 1.3pc
\hsize=6in
\input mssymb  
\vskip.4in
\def\res{\restriction}   
\def\A{\hbox{$\cal{A}$}}         
         
\def\G{\hbox{$\cal{G}$}}        \def\II{\hbox{$\cal{I}$}} 
         
\def\L{\hbox{$\cal{L}$}}                
                
\def\P{\hbox{$\cal{P}$}}                
        \def\T{\hbox{$\cal{T}$}}        
\def\NN{\hbox{$\Bbb{N}$}}

\def\al{\alpha}         \def\be{\beta}          \def\ga{\gamma}
\def\de{\delta}                
                     
\def\la{\lambda}                 
\def\vf{\varphi}        \def\om{\omega}
\def\Ga{\Gamma}                  
 
\def\bb #1{{\bar #1}}   \def\bba{\bar A}        \def\bbb{\bar B}        
\def\bbp{\bar P}        \def\bbq{\bar Q}        \def\bbu{\bar U}
\def\bbv{\bar V}        \def\bbw{\bar W}        \def\bbx{\bar X}
\def\bby{\bar Y}        \def\bbz{\bar Z}		\def\bbh{\bar P}
\def\vvp{\bar P^*}\def\vvu{\bar U^*}\def\vvv{\bar V^*}
\def\vvw{\bar W^*}\def\vvq{\bar Q^*}

\def\ubb{\bar U^@}\def\wbb{\bar W^@}\def\qbb{\bar Q^@}
\def\pbb{\bar P^@}
\def\qmu{\hbox{$Q$\llap{\lower5pt\hbox{$\sim$}\kern1pt}$_\mu$}}
\def\mo{\models}        \def\fo{\ ||\kern -4pt-}                
\def\imp{\Rightarrow}   
\def\sb{\subseteq}      \def\sbb{\subset} 
\def\em{\emptyset}      \def\sm{\setminus}
\def\lan{\langle}       \def\ran{\rangle}
\def\all{\forall}       \def\ex{\exists}
\def\and{\ \&\ }        \def\ldl{,\ldots,} 
\def\ee{\thinspace^\wedge} 
\def\l{\ell}            \def\ale{\aleph_0} 

\def\th{\hbox{\rm Th}}    \def\ath{\hbox{\rm $a$-Th}}
\def\lath{\hbox{\rm $(\la\sm a)$-Th}}
\def\cf{{\rm cf}}               \def\lg{{\rm lg}}                                              
\def\dpp{{\rm dp}}              
        
\def\max{\hbox{\rm max}}        \def\min{\hbox{\rm min}}
                  
\def\\{\noindent}               \def\sec{\par \ \par}
\def\vv{\par\rm}                \def\vvv{\medbreak\rm}

\def\pt #1. {\medbreak{\bf#1. }\enspace\sl}
\def\pd #1. {\medbreak{\bf#1. }\enspace}
\def\dd{Definition\ }  
      
\def\qed{\line{\hfill$\heartsuit$}\medbreak}
\def\sc{semi--club}                     \def\st{such that\ } 
\def\sg{semi-homogeneous\ }             
\def\today{\ifcase\month\or January\or February\or March\or April\or May
\or June\or July\or August\or September\or October\or November\or December\fi
\space\number\day, \number\year}
\def\sh{[Sh]}

\def\gu{[Gu]}
\def\ls{[LiSh]}

\font\bfa=cmbx12 scaled\magstep1        
\font\title=cmbx12 
     \let\he=\heads
    \let\hee=\headss 
\def\hh #1\par{\underbar{\he#1}\vskip1pc}
\def\hhh #1\par{\underbar{\hee#1}\vskip1.2pc} 
\def\tit #1\par{\centerline{\title #1}} 
\def\ttit #1\par{\centerline{\bfa #1}} 
\font\author=cmr10
\font\adress=cmsl10 
\def\Bby{\centerline{\by BY}\par}
\def\aut #1\par{\centerline{\author #1}}

\font\abstract=cmr8 
\font\by=cmr8           
        
\def\ar #1\ {\bf#1}

\def\abs{\abstract\centerline{ABSTRACT}\par}
\def\abss #1\par{\abstract\midinsert\narrower\narrower\noindent #1\endinsert}
\font\ba=cmr8
\font\bs=cmbxti10
\font\bib=cmtt12

\tit Random Graphs in the Monadic Theory of Order.
\par
\vskip.5in
\Bby
\aut SHMUEL LIFSCHES and SAHARON SHELAH\footnote*{The second author would 
like to thank the U.S.--Israel Binational Science Foundation for 
partially supporting this research. Publ. 527}
\par
\centerline{\adress Institute of Mathematics, The Hebrew University, 
Jerusalem, Israel} \par
\vskip 1in
\abs
\abss We continue the works of Gurevich-Shelah and Lifsches-Shelah by showing 
that it is consistent with ZFC that the first-order theory of random graphs 
is not interpretable in the monadic theory of all chains. It is provable from 
ZFC that the theory of random graphs is not interpretable in the monadic 
second order theory of short chains (hence, in the monadic theory of the real 
line).  \par
\vskip1in
\hh 0. Introduction \par \rm
We are interested in the monadic theory of 
order -- the collection of monadic sentences that are satisfied by 
every chain (= linearly ordered set). 
The monadic second-order logic is the fragment of the full second-order
logic that allows quantification over elements and over monadic (unary)
predicates only. The monadic version of a first-order language $L$ can be
described as the augmentation of\/ $L$ by a list of quantifiable set 
variables and by new atomic formulas $t\in X$ where $t$ is a first order 
term and $X$ is a set variable. 

It is known that the monadic theory of order and the monadic theory of the 
real line are at least as complicated as second order logic ([GuSh2], [Sh1]). 
The question that we are dealing with in this paper is related to the 
expressive power of this theory: what can be interpreted in it? 

In our notion of (semantic) interpretation, interpreting a theory $T$ 
in the monadic theory of order is defining models of\/ $T$ in chains. Some 
problems about the interpretability power of the monadic theory of order, 
which is a stronger criterion for complicatedness, have been raised and 
answered. For example, second order logic was shown to be even interpretable in the 
monadic theory of order ([GuSh3]) but this was done by using a weaker, 
non-standard form of interpretation: into a Boolean valued model. 

Using standard interpretation ([GMS]) it was shown that it is consistent that 
the second--order theory of\/ $\om_2$ is interpretable in the monadic theory 
of\/ $\om_2$ (hence in the monadic theory of well orders). 
On the other hand, by [GuSh], Peano arithmetic is not interpretable 
in the monadic theory of short chains, (chains that do not embed $(\om_1,<)$  
and $(\om_1,>)$) and in particular in the monadic theory of the real line. 
In [LiSh] we filled the gap left by the previous results and showed that it 
is  not provable from ZFC that Peano arithmetic is interpretable in the 
monadic theory of order. 

Here we replace Peano arithmetic by a much simpler theory -- the theory of 
random graphs, and obtain the same results by proving: 

\pt Theorem. There is a forcing notion $P$ such that in $V^P$, the theory  
of random graphs is not interpretable in the monadic second-order theory of 
chains. 
\vvv 

In fact we show that the model $V^P$ in which Peano arithmetic is not 
interpreted is a model in which the theory of random graphs is not 
interpreted (an exact formulation of the non-interpretability theorem 
is given in section 2). 

The proof is similar in its structure to the proof in [LiSh]: 
we start by defining, following [Sh], our basic objects of manipulation -
partial theories. Next, we present the notion of interpretation and  
the main theorem. We show in $\S3$ that an interpretation in a 
chain $C$ `concentrates' on an initial segment $D\sb C$ called a major 
segment. One of the main differences from [GuSh] and [LiSh] is that the 
notion of a major segment is not as sharp as there; this results in 
the need to apply more complicated combinatorial arguments. 

The most widely used idea in the proof is applying the operation of shuffling 
subsets $X,Y\sb C$: given a partition of\/ $C$, \ $\lan S_j: j\in J \ran$ and 
a subset $a\sb J$, the shuffling of\/ $X$ and $Y$ with respect to $J$ and $a$ 
is the set: \ $\bigcup_{j\in a}(X\cap S_j)\cup\bigcup_{j\not\in a}(Y\cap S_j)$. 
One of the main results in [LiSh] was to show that this operation 
preserves partial theories; this is stated and used here as well. 

To prove the main theorem we try to derive a contradiction from the existence 
of an interpretation in a chain $(C,<)\in V^P$. We start by making two  
special assumptions: that $C$ itself is the minimal major initial segment, 
and that $C$ is an uncountable regular cardinal. The spirit of the proof 
and main tools are similar to [LiSh], but some of the techniques have to be
more tortuous. The proof in this case contains all the main ingredients and
disposing of the special assumptions is essentially a formality.

Although we use many definitions and techniques from [Sh], [GuSh] and [LiSh] 
we have tried to make this paper as self contained as possible. 
The only main proof we have omitted is that of the theorem on 
preservation of partial theories under shufflings, as its proof is quite long 
and involves ideas that are not directly related to this paper. 
\sec\sec
\hh 1. Composition and Preservation of Partial Theories                 \par 
In this section we define formally the monadic theory of a chain and our main    
objects of interest: its finite approximations (partial theories). We state 
the useful properties of partial theories, namely the composition theorem and 
the theorem about preservation under shuffling. 

\

The monadic theory of a chain is defined to be the first order theory of its 
power set.

\pd \dd 1.1. Let $(C,<)$ be a chain. The {\sl monadic second-order theory 
of\/ $C$} is the first-order theory of the model 
$$C^{\rm mon}=(\P(C)\,;\;\sb,\,<^*,\,{\rm EM},\,{\rm SING})$$
where $\P(C)$ is the power set of\/ $C$, $<$ and $\sb$ are binary relations,
SING and EM are unary relations and:

$(i)$ \ $C^{\rm mon}\mo\,$SING$(X)$ iff\/ $X$ is a singleton,

$(ii)$ \ $C^{\rm mon}\mo X<^*Y$ iff\/ $X=\{x\}$, $Y=\{y\}$ (where $x,y\in C$) 
and $C\mo x<y$,

$(iii)$ \ $C^{\rm mon}\mo$EM$(X)$ iff\/ $X=\em$,

$(iv)$ \ $\sb$ is interpreted as the usual inclusion relation between 
subsets of\/ $C$.
\vvv

\pd Remark. We denote the first order language above by $L(mon)$. However we 
will be slightly informal about that and identify it with the monadic version 
of the first-order language of order, $\L$. 
Now each $\vf\in\L$ can be translated to a first-order formula 
$\vf'\in L(mon)$ by the rules: $(\ex x)\psi(x)$ (individual quantification) 
will be translated to $(\ex X)[{\rm SING}(X)\,\&\,\psi'(X)]$ and $x\in Y$ to 
${\rm SING}(X)\,\&\,(X\sb Y)$. So when we write $C\mo\vf$ (for $\vf\in\L$) we 
mean $C^{\rm mon}\mo\vf'$ and $x<y$ is translated as $X<^*Y$. 

\

\pd Notations 1.2. We denote individual variables by $x,y,z$ and set variables 
by $X,Y,Z$. $a,b,c$ are elements and $A,B,C$ are sets. $\bb a$ and $\bba$ 
denote finite sequences having lengths $\lg(\bb a)$ and $\lg(\bba)$.  
We will write $\bb a\in C$ and $\bba\sb C$ instead of 
$\bb a\in\,^{\lg(\bb a)}C$ or $\bba\in\,^{\lg(\bba)}\P(C)$, we may also write   
$a_0\in\bb a$ or $A_0\in\bba$.
\vv

\

Next is the definition of the partial $n$-theory of\/ $\bba$ in $C$

\pd \dd 1.3. Let $(C,<)$ be a chain and $\bba\sb C$. We define 
$$t = \th^n(C;\bba)$$ 
by induction on $n$:

\\for $n=0$:\ \ \ \ \ \ \ \ \ \ \ $t=\{\vf(\bbx):\vf\in L(mon), \ \vf {\rm \ 
quantifier \ free}, \ C^{\rm mon}\mo\vf(\bba)\}$ 

\\for $n=m+1$:
\ \ \ $t = \{\th^m(C;\bb A\ee B):B\sb C\}$. 
\vvv

\

\pt Lemma 1.4. (A) For every formula $\psi(\bbx)\in\L$ there is an $n$  
such that from $\th^n(C;\bba)$ we can decide effectively whether 
$C\mo\psi(\bba)$. We call the minimal such $n$ {\sl the depth of\/ $\psi$} 
and write $\dpp(\psi)=n$.

\\(B) \ For every $n$ and $l$ there is a finite set of monadic formulas 
(effectively computable from $n$ and $l$) $\Psi(n,l)=
\{\psi_m(\bbx):m<m^*,\; \lg(\bbx)=l\}\sb\L$ such that for any 
chains $C$, $D$ and $\bba\sb C$, $\bbb\sb D$ of length $l$ the 
following hold: 

(1) \ $\dpp(\psi_m(\bbx))\le n$ for $m<m^*$,

(2) \ $\th^n(C;\bba)$ can be computed from $\{m<m^*:C\mo\psi_m[\bba]\}$, 

(3) \ $\th^n(C;\bba)=\th^n(D;\bbb)$ iff for every
    $m<m^*$, \ \ $C\mo\psi_m[\bba]\iff D\mo\psi_m[\bb B]$. 
\vv

\pd Proof. In \sh, \ Lemma 2.1.
\vv \qed

\

\pd \dd 1.5. When $\Psi(n,l)$ is as in 1.4(B), for each chain $C$ 
and $\bba\sb C$ of length $l$ we can identify $\th^n(C;\bba)$
with a subset of\/ $\Psi(n,l)$. Denote by $T_{n,l}$ the collection of
subsets of\/ $\Psi(n,l)$ that arise from some $\th^n(C;\bba)$ 
and call it {\sl the set of formally possible $(n,l)$-theories}. 
\vvv

\

\pd Remark. For given $n,l\in\NN$, each $\th^n(C;\bba)$  
is hereditarily finite, (where $\lg(\bba)=l, \ C$\ is a chain), \ and we 
can effectively compute the set of formally possible theories $T_{n,l}$. 
(See \sh, \ Lemma 2.2). 
\vvv   

\

\pd \dd 1.6. If\/ $(C,<_C)$ and $(D,<_D)$ are chains then $(C+D,<)$ is the 
chain that is obtained by adding a copy of\/ $D$ after $C$ (where $<$ is 
naturally defined).

If\/ $(I,<)$ is a chain and $\lan(C_i,<_i):i\in I\ran$ is a sequence of 
chains then $\sum_{i\in I}(C_i,<_i)$ \ is the chain that is the 
concatenation of the $C_i$'s along $I$ equipped with the obvious order.
\vvv

\

Given $\bba=\lan A_0\ldl A_{l-1}\ran$ and 
$\bbb=\lan B_0\ldl B_{l-1}\ran$ we denote by $\bba\cup\bbb$ the sequence 
$\lan A_0\cup B_0\ldl A_{l-1}\cup B_{l-1}\ran$.
The heavily used composition theorem for chains states that 
the partial theory of a chain is determined by the partial theories of 
its convex parts. 

\pt Theorem 1.7 (Composition theorem for chains). \par
(1) If \ $C$, $C'$, $D$ and $D'$ are chains, $\bba\sb C$, $\bba'\sb C'$, 
$\bbb\sb D$ and $\bbb'\sb D'$ are of the same length and if

$$\th^m(C;\bba) = \th^m(C';\bba')$$ 
and 
$$\th^m(D;\bbb) = \th^m(D';\bbb')$$
then 
$$\th^m(C+D;\bba\cup\bbb) = \th^m(C'+D';\bba'\cup\bbb').$$ 
(2) If\/ $I$ is a chain and \ $\th^m(C_i;\bba^i)=\th^m(D_i;\bbb^i)$ for each 
$i\in I$ (with all sequences of subsets having the same length) then
$$\th^m\Big(\sum_{i\in I}C_i;\cup_i\bba^i\Big)= 
\th^m\Big(\sum_{i\in I}D_i;\cup_i\bbb^i\Big).$$ 
\vv

\pd Proof. By \sh\ Theorem 2.4 (where a more general theorem is proved), 
or directly by induction on $m$. See also theorem 1.9 below.
\vv \qed

\

Using the composition theorem we can define a formal operation of addition 
of partial theories.

\pd Notation 1.8. 

\\(1) When $t_1,t_2,t_3\in T_{m,l}$ for some $m,l\in\NN$, then 
$t_1+t_2=t_3$ means: there are chains $C$ and $D$, and $\bba\sb C$, 
$\bb B\sb D$  \st

\centerline{$t_1=\th^m(C;A_0\ldl A_{l-1}) \and t_2=\th^m(D;B_0\ldl B_{l-1}) 
\and t_3=\th^m(C+D;\bba\cup\bbb).$}

\\(By the composition theorem, the choice of\/ $C$ and $D$ is immaterial).

\\(2) $\sum_{i\in I}\th^m(C_i;\bba^i)$ \ is \ 
$\th^m\big(\sum_{i\in I}C_i;\cup_{i\in I}\bba^i)$, 
(assuming $\lg(\bba^i)=\lg(\bba^j)$ for $i,j\in I$). 

\\(3) If\/ $D$ is a sub-chain of\/ $C$ and $\bba\sb C$ then 
$\th^m(D;\lan A_0\cap D,A_1\cap D,\ldots\ran)$ is abbreviated by
$\th^m(D;\bba)$. 

\\(4) For $a<b\in C$ and $\bbp\sb C$ we denote by 
$\th^n(C;\bb P)\res_{[a,b)}$ the theory $\th^n([a,b);\bb P\cap[a,b))$.
\vvv 

\

We conclude this part by giving a monadic version of the Feferman-Vaught 
theorem. Note that the composition theorem is a consequence. 

\pt Theorem 1.9. For every $n,l<\om$ there is $m=m(n,l)<\om$, effectively 
computable from $n$ and $l$, \st if 

$(i)$ \/ $I$ is a chain,

$(ii)$ \/ $\lan C_i:i\in I\ran$ \/ is a sequence of chains,

$(iii)$ for $i\in I$, $\bbq_i\sb C_i$ is of length $l$,

$(iv)$ \/ for $t\in T_{n,l}$, \/ $P_t:=\{i\in I: \th^n(C_i;\bb Q_i)=t\}$, 

$(v)$ \/ $\bbp:=\lan P_t: t\in T_{n,l}\ran$,

\\then $\th^n(\sum_{i\in I}C_i;\cup\bb Q_i)$ is computable from 
\/ $\th^m(I;\bbp)$.
\vv

\pd Proof. This is theorem 2.4. in [Sh].
\vv\qed

\

\\Next we define semi--clubs and shufflings and we quote the important 
preservation theorem.

\pd \dd 1.10. Let $\la>\ale$ be a regular cardinal

\\1) We say that $a\sb\la$ is a 
{\sl \sc\ subset of}\/ $\la$ if for every $\al<\la$ with $\cf(\al)>\ale$:  \par
\\if\/ $\al\in a$ then there is a club subset of\/ $\al,\ C_\al$ \st $C_\al\sb a$ 
\ and \par
\\if\/ $\al\not\in a$ then there is a club subset of 
$\al,\ C_\al$ \st $C_\al\cap a=\emptyset$.  

\\(Note that $\la$ and $\emptyset$ are semi--clubs
and that a club $J\sb\la$ is a \sc\ provided that the first and the
successor points of\/ $J$ are of cofinality $\leq\ale$. Also, if\/ $a\sb\la$ is
a \sc\ then $\la\sm a$ is one as well.)

\\2) Let $X,Y\sb\la, \ J =\{\al_i:i<\la\}$ a club subset of\/ $\la$, and let
$a\sb\la$ be a \sc\ of\/ $\la$. We will define the {\sl shuffling of
$X$ and $Y$ with respect to $a$ and $J$},\ denoted by $[X,Y]^{J}_a$, as:
  $$[X,Y]^{J}_a = \bigcup_{i\in a}\big(X\cap[\al_i,\al_{i+1})\big)\cup
                     \bigcup_{i\not\in a}\big(Y\cap[\al_i,\al_{i+1})\big)$$

\\3) When $\bb X,\bb Y\sb\la$ are of the same length, we define 
$[\bb X,\bb Y]_a^J$ naturally.                                    

\\4) We can naturally define shufflings of subsets of an ordinal $\de$ with 
respect to a club $J\sb\de$ and a \sc\  $a\sb{\rm otp}(J)$.     

\\6) $\ath^n(\la;\bb P)$ is $\th^n(\la;\bb P,a)$ where $a\sb\la$ is a \sc. 
\vvv

\

The next theorem, which will play a crucial role in contradicting the 
existence of interpretations, states that the result of the shuffling of 
subsets of the same type is an element with the same partial theory. The 
proof of the preservation theorem in $\S$4 of \ls\ requires some amount of 
computations and uses some auxiliary definitions that are not material in 
the other parts of the paper. For example, the partial theories 
WTh$^n(C;\bb P)$, ATh$^n(\be,\ (C;\bb P))$ and $a$-WA$^n(C;\bb P)$ are used 
in the proof and even the formulation of the theorem but we can avoid 
defining them  by noticing that (for a large enough $m$) $\ath^m(C;\bb P)$ 
computes all these partial theories. We also avoid the definition of an 
$n$-suitable club (which relies on ATh$^n$). All the details can be found 
of course in \ls.

\pt Theorem 1.11 (preservation theorem). Let $\bb P_0,\bb P_1\sb\la$ be of 
length $l$, $n<\om$ and $a\sb\la$ be a \sc.  \par
\\Then there are an $m=m(n,l)<\om$ and a club $J=J(n,\bb P_0,\bb P_1)\sb\la$ 
\st if\/ $\bb X := [\bb P_0,\bb P_1]^J_a$ then 
$$(*) \ \Big[\ath^m(\la;\bb P_0) = \ath^m(\la,\bb P_1)\Big] \imp 
\Big[\th^n(\la;\bb P_0) = \th^n(\la;\bb P_1) = \th^n(\la;\bb X)\Big].$$
Moreover, there is $t^*=t^*(\bb P_0,\bb P_1)\in T_{n,l}$ such that,
for every $\ga\in J$ with $\cf(\ga)=\ale$, 
$$(**)\hskip3cm\Big[\ath^m(\la;\bb P_0)=\ath^m(\la;\bb P_1)\Big] 
\imp\hskip3cm$$
$$\Big[\th^n(\la;\bb P_0)\res_{[0,\ga)}=\th^n(\la;\bb P_1)\res_{[0,\ga)}=
\th^n(\la;\bb X)\res_{[0,\ga)} = t^*\Big].$$ \vv
\pd Proof. By \ls\ 4.5, 4.12. 
\vv \qed

\

\pd \dd 1.12. Let $\bb P_0,\bb P_1\sb\la$ be as above. Call a club $J\sb\la$ 
{\sl an $n$-suitable club for $\bb P_0$ and $\bb P_1$} if for every \sc\  
$a\sb\la$, $(*)$ and $(**)$ of 1.11 hold. 
\vv 

\

\pt Fact 1.13. For every finite sequence 
{\bf P}=$\lan \bb P_i:\bb P_i\sb\la,\ \lg(\bb P_i)=l,\ i<k\ran$ 
and for every $n<\om$ there is a club $J\sb\la$ that is $n$-suitable for 
every pair from {\bf P}. \vv
\pd Proof. By \ls\ 4.3, 4.4.
\vv \qed
\sec\sec
\hh 2. Random Graphs and Uniform Interpretations \par

The notion of semantic interpretation of a theory $T$ in a theory $T'$ is not 
uniform. Usually it means that models of\/ $T$ are defined inside models of 
$T'$ but the definitions vary with context. In \ls\ we gave the general 
definition of the notion of interpretation of one first order theory in 
another. In our case, in which we deal with interpreting a class of theories,   
another notion emerges, that of a {\sl uniform interpretation}. 

\

First we define the theory of\/ $K$-random graphs:

\pd \dd 2.1. Let $1<K\le\om$. An undirected graph $\G=(G,R)$ is a 
{\sl $K$-random graph} if 
$$\Bigl[A_0,A_1\sb G\ \&\ |A_0|,|A_1|<K\ \&\ A_0\cap A_1=\em\Bigr] \ \imp \ 
\Bigl[(\ex x\in G)(\all a\in A_0)(\all b\in A_1)[xRa\ \&\ \neg xRb]\Bigr].$$
(When this holds we will say that $x$ {\sl separates} $A_0$ from $A_1$). 
\vvv

\

\pd \dd 2.2. (1) $RG_K$ is the theory of all $K$-random graphs (that is  
all the sentences, in the first-order language of graphs, that are 
satisfied by every $K$-random graph). 
$RG_K^i$ is theory of all the infinite graphs that are $K$-random. 

\

\\(2) $\Ga_K$ is the class of all the $K$-random graphs, (clearly 
$K<L\le\om\imp\Ga_L\sb\Ga_K$). $\Ga_K^i$ is the classes of infinite 
$K$-random graphs.

\

\\(3) $\Ga_{\rm fin}$ is the class $\{\Ga_K\}_{1,<K<\om}$,  
$\Ga_{\rm fin}^i$ is $\{\Ga_K^i\}_{1,<K<\om}$.
\vvv

\

The next definition is the one used in \ls.  It is applicable in dealing with 
$RG_\om$, but will have to be modified for dealing with finitely-random 
graphs.
 
\pd \dd 2.3.  An {\sl interpretation} of a model $\G$\ of\/ $RG_K$ in the 
monadic theory of a chain $C$ is a sequence of formulas in the language $\L$ 
of the monadic theory of order                                      
$$\II = \lan U(\bb X,\bb W), \ E(\bb X,\bb Y,\bb W), \ 
R(\bb X,\bb Y,\bb W) \ran$$
\\where: 

1) $U(\bb X,\bb W)$ is the {\sl universe formula} that says which sequences 
of subsets of\/ $C$ represent elements of\/ $\G$.  
We denote by $C^U$ the set $\{\bb X\sb C : C\mo U(\bb X,\bb W)\}$.    

2) $E(\bb X,\bb Y,\bb W)$ is the {\sl equality formula}, an equivalence 
relation on $C^U$. 
We write $\bb A\sim\bb B$ when $C\mo E(\bb A,\bb B,\bb W)$. 

3) $R(\bb X,\bb Y,\bb W)$ is the {\sl interpretation of the graph relation}, 
a binary relation on $C^U$ which respects $\sim$ i.e. 
``$C\mo R(\bb A,\bb B,\bb W)$''  depends only on the $E$-equivalence 
classes of\/ $\bb A$ and $\bb B$.

4) $\bb W\sb C$ is a finite set of parameters allowed in the interpreting  
formulas. 

5) $\lan C^U/\sim\,, \ R \ran$ \ $\cong$ \ $\G$. 

\vvv

\

\pd \dd 2.4. Let $\II$ be an interpretation of\/ $\G$ in the 
monadic theory of a chain $C$.  

The {\sl dimension} of the interpretation, denoted by $d(\II)$,   
is $\lg(\bbx)$. We will usually assume without loss of generality 
that $\lg(\bbw)=d(\II)$ as well. 

The {\sl depth} of the interpretation, denoted by $n(\II)$, is 
$\max\{\dpp(U),\dpp(E),\dpp(R)\}$. 
\vvv

\

\pd \dd 2.5. Let $RG^*$ be one of the theories defined in 2.2(1) and $\Ga^*$ 
be the respective class. We say that {\sl the monadic theory of order 
interprets $RG^*$ (or $\Ga^*$)} if there is a chain $C$, a random graph 
$\G\in\Ga^*$ and an interpretation $\II=\lan U(\bbx,\bbw), \ E(\bbx,\bby,\bbw),
 \ R(\bbx,\bby,\bbw)\ran$ with $\lan C^U/\sim, \ R \ran$ \ $\cong$ \ $\G$. 
\vvv

\

Common notions of an interpretation of a theory $T_1$ in a theory $T_2$ 
demand that every model of\/ $T_1$ is interpretable in a model of\/ $T_2$ (as in 
[BaSh]) or that inside every model of\/ $T_2$ there is a definable model of 
$T_1$ (see [TMR]). Here we seem to require the minimum: a single model is 
interpreted in a single chain. This is often useful, but not always:

\pt Fact 2.6. For every $1<K<\om$ there is a chain $C$ and a sequence $\II$ 
\st $\II$ is an interpretation of a model of\/ $RG_K^i$ (hence of\/ $RG_K$) in 
$C$. (That is, $RG_K$, $RG_K^i$, $RG_{\rm fin}$, $RG^i_{\rm fin}$  are 
interpretable in the monadic theory of order). 
\vv

\pd Proof. We shall demonstrate the construction for $K=2$; the other cases 
are similar. 

Let $C=(\om,<)$, we will show that there is a one-dimensional interpretation 
of an infinite model of\/ $RG_2$ in $C$ without parameters. For that we have 
to define $U(X)$, $E(X,Y)$ and $R(X,Y)$. Let:

\\$U(X)$ := $[X=\{a\}\ \&\ a>1] \vee [X=\{x,a,b\}\ \&\ x\in\{0,1\}\ \&\ 
a,b>1] \ ;$

\\$E(X,Y)$ := $U(X)\ \&\ U(Y)\ \&\ X=Y \ ;$

\\$R(X,Y)$ := $U(X)\ \&\ U(Y)$ and either: 

$[X=\{a\}\ \&\ Y=\{0,a,b\}\ \&\ a<b]$ or 

$[Y=\{a\}\ \&\ X=\{0,a,b\}\ \&\ a<b]$ or 

$[X=\{b\}\ \&\ Y=\{1,a,b\}\ \&\ a>b]$ or 

$[Y=\{b\}\ \&\ X=\{1,a,b\}\ \&\ a>b]$ or

$[X=\{a\}\ \&\ Y=\{x,c,d\}\ \&\ x\in \{0,1\}\ \&\ a\not\in\{c,d\}]$ or

$[Y=\{a\}\ \&\ X=\{x,c,d\}\ \&\ x\in \{0,1\}\ \&\ a\not\in\{c,d\}]$.

\\Clearly everything is expressible in $\L$ and $R(X,Y)$ defines on 
$\{X\sb\om: (\om,<)\mo U(X)\}$ a graph relation that is 2-random.
\vv \qed

\

Motivated by the previous fact we will define now the suitable modification of 
the previous definitions. The idea is to interpret, in a uniform way, 
an infinite set of random graphs.  

\pd \dd 2.7.  A {\sl uniform interpretation} of\/ $\Ga_{\rm fin}$ in the monadic 
theory of order is a sequence 
$$\{\lan\,C_K,\ I,\ \bbw_K\,\ran : K\in A\}$$ 
\\where \par
\\1) $C_K$ is a chain, \par
\\2) $A$ is an infinite subset of\/ $\om$, \par
\\3) $\bb W_K\sb C_K$ for $K\in A$, \par
\\4) $I=\lan U(\bb X,\bb Z),E(\bb X,\bb Y,\bb Z),R(\bb X,\bb Y,\bb Z)\ran$
is a sequence of formulas in $\L$, \par 
\\5) $\II_K:=\lan U(\bbx,\bbw_K),E(\bbx,\bby,\bbw_K),R(\bbx,\bby,\bbw_K)\ran$
is an interpretation of a model of\/ $RG_K$ in $C_K$ for $K\in A$.  
\vvv

\

Given $d$ and $n$ in $\NN$ there is only a finite number of possible 
interpretations $\II$ having dimension $d$ and depth $n$. The following is 
therefore clear:

\pt Proposition 2.8. The following are equivalent:

\\{\rm(A)} There is no uniform interpretation of\/ $\Ga_{\rm fin}$ in the 
monadic theory of order. 

\\{\rm(B)} For every $n,d\in\NN$ there is $K^*=K^*(n,d)\in\NN$ \st 
if\/ $K\ge K^*$ and $\II$ is an interpretation of some $\G\mo RG_K$ in a 
chain $C$, then either $d(\II)>d$ or $n(\II)>n$.

\\{\rm(C)} For every sequence $\II=\lan U(\bbx,\bbz), 
\ E(\bbx,\bby,\bbz), \ R(\bbx,\bby,\bbz)\ran$ there is $K^*=K^*(\II)<\om$ \st 
there are no chain $C$, $\bb W\sb C$, $K\ge K^*$ and $\G\in\Ga_K$ \st 
$\lan U(\bbx,\bbw),E(\bbx,\bby,\bbw),R(\bbx,\bby,\bbw\ran$
is an interpretation of\/ $\G$ in $C$. 
\vv\qed

\

Our main theorem has therefore the following form:

\pt Theorem 2.9 (Non-Interpretability Theorem). There is a forcing notion $P$ 
such that in $V^P$ the following hold:  

\\(1) $RG_\om$ is not interpretable in the monadic theory of order.

\\(2) For every sequence of formulas $\II=\lan U(\bbx,\bbz), 
\ E(\bbx,\bby,\bbz), \ R(\bbx,\bby,\bbz)\ran$ there is $K^*<\om$, 
(effectively computable from $\II$), \st for no chain $C$, $\bb W\sb C$, and 
$K\ge K^*$ does $\lan U(\bbx,\bbw),$ $E(\bbx,\bby,\bbw),R(\bbx,\bby,\bbw)\ran$
interpret $RG_K$ in $C$. 

\\(3) The above propositions are provable in ZFC. if we restrict ourselves 
to the class of short chains. 
\vv

\

\pd Remark. As an $\om$-random graph is $K$-random for every $K<\om$, an 
interpretation of\/ $RG_\om$ is a uniform interpretation of\/ $\Ga_{\rm fin}$.   
Therefore clause (1) in the non-interpretability theorem follows from 
clause (2).
\vv
\sec\sec 
\hh 3. Major and Minor Segments                                    \par 
>From now on we will assume that there exists (in the generic model $V^P$ that 
is defined later) a uniform interpretation $\II$\/ of\/ $\Ga_{\rm fin}$ in the 
monadic theory of order. For reaching a contradiction we have to find a large 
enough $K=K(\II)<\om$ (a function of the depth and dimension of\/ $\II$) and 
show that no chain interprets a $K$-random graph by $\II$. 
The aim of this (and the next) section is to gather facts that will enable us 
to compute an appropriate $K$.  The main observation is that an interpretation 
in a chain $C$ ``concentrates'' on a segment (called a {\sl major segment}).  
One of the factors in determining the size of\/ $K$ will be the relation between 
the major segment and the other, {\sl minor}, segment.

\

\pd Context 3.1. 
$\II=\lan U(\bb X,\bb W),\ E(\bb X,\bb Y,\bb W),\ R(\bb X,\bb Y,\bb W) \ran$  
is an interpretation of a $K$-random graph $\G=(G,R)$ on a chain $C$. 
$\bb W\sb C$ are the parameters, $d=d(\II)=\lg(\bbx)=\lg(\bbw)$ is the 
dimension of\/ $\II$ and $n=n(\II)$ is its depth. 
\vvv

\

\pd \dd 3.2. $A\sb G$ is {\sl big for} $(K_1,K_2)$ if there is $B\sb G$ 
with $|B|\le K_1$ \st:  \par 
\\$(*)$ for every disjoint pair $A_1,A_2\sb G\sm B$ with 
$|A_1\cup A_2|\le K_2$ there is some $x\in A\sm(A_1\cup A_2)$ that separates 
$A_1$ from $A_2$ i.e.   
$\bigl(\bigwedge_{y\in A_1}xRy\bigr)\wedge
\bigl(\bigwedge_{y\in A_2}\neg xRy\bigr)$. 

\\When $(*)$ holds we say that $B$ witnesses the $(K_1,K_2)$-bigness of\/ $A$.
\vv

\

Non-bigness is an additive property: 

\pt Proposition 3.3. Let $A\sb G$ be big for $(K_1,K_2)$ and suppose that
$A=\bigcup_{i<m}A_i$. Then there is an $i<m$ \st $A_i$ is big for 
$(K_1+K_2,K_2/m)$. 
\vv
\pd Proof. Let $B\sb G$ ($|B|\le K_1$) witness the bigness of\/ $A$. For $i<m$ 
we will try to define by induction counter-examples for bigness, that is a set 
$B_i\sb G$ and a function $h_i$ so that: 

(1) $|B_i|\le K_2/m$,      \par
(2) $B_i\sb G\sm(B\cup\bigcup_{j<i}B_j)$,    \par
(3) $h_i\colon B_i\to\{t,f\}$,     \par
(4) for no $x\in A_i\sm B_i$ we have 
$(\all y\in B_i)[xRy\leftrightarrow h(y)=t]$.

\\Suppose we succeed. Let $C_1:=\big\{x: \bigvee_i(x\in B_i\ \&\ h_i(x)=t) \big\}$ 
and $C_2:=\big\{x: \bigvee_i(x\in B_i\ \&\ h_i(x)=f) \big\}$. 
But $|C_1\cup C_2|\le K_2$,\ $C_1\cup C_2\sb G\sm B$ and of course \
$C_1\cap C_2=\em$ so by the assumption on $A$ there is some $x\in A$ that 
separates $C_1$ from $C_2$. Such an $x$ belongs to some $A_i$ and it 
separates $C_1\cap B_i$ from $C_2\cap B_i$. This contradicts clause (4).

Therefore at some stage $i<m$ we can't define $B_i$ and look at 
$B^*:=B\cup\bigcup_{j<i}B_j$. Now $|B^*|\le K_1+i\cdot K_2/m\le K_1+K_2$ and  
``being unable to continue'' means: if\/ $B_i\sb G\sm B^*$ and 
$|B_i|\le K_2/m$ then for every partition of\/ $B_i$ to $B_i^1$ and $B_i^2$ 
there is some $x\in A_i\sm B_i$ such that 
$\bigl(\bigwedge_{y\in {B_i^1}}xRy\bigr)\wedge
\bigl(\bigwedge_{y\in {B_i^2}}\neg xRy\bigr)$.
In other words, $A_i$ is big for $(K_1+K_2,K_2/m)$ (witnessed by $B^*$) as 
required.
\vv \qed

\

\pd Notation 3.4. $\bb A\sb C$ is called {\sl a representative} if it 
represents an element of\/ $\G$ i.e. if\/ $C\mo U(\bb A,\bb W)$ 
(of course $\lg(\bb A)=d)$. The representatives $\bb A,\bb B\sb C$ are called 
{\sl equivalent} and we write $\bb A\sim\bb B$ if they represent the same 
element in $\G$ i.e. if\/ $C\mo E(\bb A,\bb B,\bb W)$. We use upper case letters 
such as $\bb X, \bb A, \bb U_i$ to denote representatives. The corresponding 
lower case letters ($x, a, u_i$) will denote the elements of\/ $\G$ that are 
represented by the former. So e.g. $\bb A\sim \bb B\iff a=b$.
\vvv

\pd \dd 3.5. \par
\\1) A sub-chain $D\sb C$ is a {\sl segment} if it is convex (i.e. 
$x<y<z \ \&\ x{,}z\in D\ \imp \ y\in D$).                        

\\2) A {\sl Dedekind cut} of\/ $C$ is a pair $(L,R)$ where $L$ is an initial 
segment of\/ $C$, \ $R$ is a final segment of\/ $C$, \ $L\cap R=\em$ and
$L\cup R=C$.                                                    

\\3) Let $\bb A,\bb B\sb C$. We will say that $\bb A,\bb B$ 
{\sl coincide on ({\rm resp.} outside)} a segment $D\sb C$, 
if\/ $\bb A\cap D=\bb B\cap D$ (resp. $\bb A\cap(C\sm D)=\bb B\cap(C\sm D) \ )$.

\\4) The {\sl bouquet size} of a segment $D\sb C$ denoted by $\#(D)$ is the 
supremum of cardinals $|S|$ where $S$ ranges over collections of 
nonequivalent representatives coinciding outside $D$. Thus $\#(D)\ge n$ iff 
there are nonequivalent representatives $A_1,A_2,\ldots,A_n$ coinciding 
outside $D$.  
\vvv

\pd \dd 3.5. Let $D\sb C$ be a segment 

\\1) $D$ is $i^*$-{\sl fat} if\/ $\#(D)\ge i^*$

\\2) $D$ is $(K_1,K_2)$-{\sl major} if there is a set $\{ \bb U_i:i<i^* \}$ 
of representatives coinciding outside $D$, and representing a subset of\/ $\G$ 
that is big for $(K_1,K_2)$.

\\3) $D$ is called $(K_1,K_2)$-{\sl minor} if it not $(K_1,K_2)$-major. 
\vvv

\

We denote by $M_1$ the number $|T_{n,3d}|$ (i.e. the number of possibilities 
for $\th^n(C;\bb X,\bb Y,\bb Z)$).

\pt Proposition 3.7. Let $(L,R)$ be a Dedekind cut of\/ $C$. If\/ $L$ [$R$] is 
$(K_1,K_2)$-major then $R$ [$L$] is not $K_3$-fat where $K_3=M_1(K_1+K_2)+1$.  
\vv

\pd Proof. Suppose $\lan \bb A_i: i<i^L\ran$ demonstrate that $L$ is 
$(K_1,K_2)$-major, i.e. they represent a $(K_1,K_2)$-big set 
$\lan a_i: i<i^L\ran$ in $\G$ and $\bb A_i\res_R=\bb A^*$. Assume towards 
a contradiction that $\lan \bb B_i: i<K_3\ran$ demonstrate that 
$R$ is $K_3$-fat (i.e. $i<j<K_3\imp b_i\ne b_j$ and $\bb B_i\res_L=\bb B^*$).
Define an equivalence relation $E_L$ on $\{0,1\ldl i^L-1\}$ by: 
$$iE_Lj\iff \th^n(L;\bb A_i,\bb B^*,\bb W)=\th^n(L;\bb A_j,\bb B^*,\bb W).$$
By the definition of\/ $M_1$, $E_L$ has at most $M_1$ equivalence 
classes.  By proposition 3.3 there is $a^L\sb\{0,1\ldl i^L-1\}$, an $E_L$ 
equivalence class, \st $\{\bb A_i: i\in a^L\}$ represents a 
$(K_1+K_2,K_2/M_1)$-big subset of $\G$. 
Let $B\sb\G$ witness the $(K_1+K_2,K_2/M_1)$-bigness 
of\/ $\{a_i: i\in a^L\}$. Since $|B|\le (K_1+K_2)$ and $K_3=M_1(K_1+K_2+1)$ 
we can choose some $j_1,j_2<K_3$ with $b_{j_1},b_{j_2}\not\in B$ and with 
$\th^n(R;\bb A^*,\bb B_{j_1},\bb W)=\th^n(R;\bb A^*,\bb B_{j_2},\bb W)$.
Now by the composition theorem 1.7, and the choice of $a^L$ and $j_1, j_2$ 
we have for every $i\in a^L$: 

$\th^n(C;\bb A_i,\bb B_{j_1},\bb W)$ = 
$\th^n(L;\bb A_i,\bb B_{j_1},\bb W)+
\th^n(R;\bb A_i,\bb B_{j_1},\bb W)$ = 

$\th^n(L;\bb A_i,\bb B^*,\bb W)$+$\th^n(R;\bb A^*,\bb B_{j_1},\bb W)$ = 
$\th^n(L;\bb A_i,\bb B^*,\bb W)$+$\th^n(R;\bb A^*,\bb B_{j_2},\bb W)$. 

\\Therefore for every $i\in a^L$ 
$$C\mo R(\bb A_i,\bb B_{j_1},\bb W)\iff R(\bb A_i,\bb B_{j_2},\bb W).$$ 
Since $b_{j_1},b_{j_2}\not\in B$ we get a contradiction to 
``$A$ is $(K_1+K_2,K_2/M_1)$-big as witnessed by $B$''.  
\vv \qed
	   
\

\pd Notation 3.8. Let $M_2$ be $|T_{n,2d}|$, $M_3$ be $M_1+1$ 
($=|T_{n,3d}|+1$) and $M_4$ be \st for every colouring 
$f\colon[M_4]^3\to\{0,1,\ldots,6\}$ there is a homogeneous subset 
of\/ $\{0,1\ldl M_4-1\}$ of size $M_3$, where $[M_4]^3$ is 
$\{\lan i,j,k\ran : \ i<j<k<M_4\}$. ($M_4$ exists by Ramsey theorem).
\vv

The main lemma states that in every Dedekind cut one segment is major. Now we 
have to make an assumption on the degree of randomness of\/ $\G$. 

\pt Lemma 3.9. Assume $K>(M_3)^2$ ($K$ is from ``$K$-random''). Let $(L,R)$ be 
a Dedekind cut of\/ $C$. Then either $L$ or $R$ is $(K_1,K_2)$-major where 
$K_1=K+{K\over{(M_2)^2}}$ \ and \ $K_2={K\over{(M_2)^2}\cdot M_4}$. \vv

\pd Proof. $\G$ is $(0,K)$-big and let $\{\bb U_i: i<i^*\}$ be a list of 
representatives for the elements of\/ $\G$. Define a pair of equivalence 
relations $E^0_L$ and $E^0_R$ on $i^*=|\G|$ by: 
$$\Big[iE^0_Lj\ \iff \ \th^n(L;\bb U_i,\bb W)=\th^n(L;\bb U_j,\bb W)\Big] \ \ \ 
\Big[iE^0_Rj\ \iff \ \th^n(R;\bb U_i,\bb W)=\th^n(R;\bb U_j,\bb W)\Big].$$
By the definition of\/ $M_2$ each relation has $\le M_2$ equivalence classes; 
therefore by 3.3 there is a subset $A_1\sb i^*$ and pair of theories 
$(t_1,t_2)$ \st $\{u_i: i\in A_1\}$ is $(K,{K\over(M_2)^2})$-big and
$$i\in A_1 \ \imp\  
\big[\th^n(L;\bb U_i,\bb W)=t_1 \ \&\ \th^n(R;\bb U_i,\bb W)=t_2\big].$$
Denote by $\bb X\ee\bb Y$ the tuple $(\bb X\res_L)\cup(\bb Y\res_R)$. 

\ 

$(\al)$ \ For $i,j\in A_1$ we have $C\mo U(\bb U_i\ee\bb U_j,\bb W)$
(hence $\bb U_i\ee\bb U_j$ is a representative).

\\Why? Because $C\mo U(\bb U_i,\bb W)$ and by the composition theorem 

\\$\th^n(C;\bb U_i,\bb W) = \th^n(L;\bb U_i,\bb W)+\th^n(R;\bb U_i,\bb W) = 
t_1+t_2 = $

\\$\th^n(L;\bb U_i,\bb W)+\th^n(R;\bb U_j,\bb W) = 
\th^n(C;\bb U_i\ee\bb U_j,\bb W)$. 

\ 

\\Define a pair of relations $E_L$ and $E_R$ on $\{\bbu_i:i\in A_1\}$ by: 
$$\bbu_iE_L\bbu_j\ \iff \ (\ex r\in A_1)\big[\bb U_i\ee\bb U_r\sim
\bb U_j\ee\bb U_r\big]$$
$$\bbu_iE_R\bbu_j\ \iff \ (\ex l\in A_1)\big[\bb U_l\ee\bb U_i\sim
\bb U_l\ee\bb U_j\big]$$

$(\be)$ \ $\bbu_iE_L\bbu_j\imp(\all r\in A_1)\big(\bbu_i\ee\bbu_r
\sim\bbu_j\ee\bbu_r\big)$ $\and$ $\bbu_iE_R\bbu_j\imp$
$(\all l\in A_1)\big(\bbu_l\ee\bbu_i\sim\bbu_l\ee\bbu_j\big)$.

\\Why? Suppose $\bbu_iE_L\bbu_j$, $\bbu_i\ee\bbu_r\sim\bbu_j\ee\bbu_r$ and let 
$r_1\in A_1$. 

\\Now $\th^n(R;\bb U_r,\bb W)$ = $t_2$ = $\th^n(R;\bb U_{r_1},\bb W)$ hence
$\th^n(R;\bb U_r,\bb U_r,\bb W)$ = $\th^n(R;\bb U_{r_1},\bb U_{r_1},\bb W)$. 
By the composition theorem 

\\$\th^n(C;\bb U_i\ee\bb U_{r_1},\bb U_j\ee\bb U_{r_1},\bb W)$ = 

\\$\th^n(L;\bb U_i,\bb U_j,\bb W)$+$\th^n(R;\bb U_{r_1},\bb U_{r_1},\bb W)$ = 
$\th^n(L;\bb U_i,\bb U_j,\bb W)$+$\th^n(R;\bb U_r,\bb U_r,\bb W)$ = 

\\$\th^n(C;\bb U_i\ee\bb U_r,\bb U_j\ee\bb U_r,\bb W)$. 

\\Therefore $\th^n(C;\bb U_i\ee\bb U_{r_1},\bb U_j\ee\bb U_{r_1},\bb W)$ = 
$\th^n(C;\bb U_i\ee\bb U_r,\bb U_j\ee\bb U_r,\bb W)$  
and hence 

\\$\bb U_i\ee\bb U_r\sim\bb U_j\ee\bb U_r\iff
\bb U_i\ee\bb U_{r_1}\sim\bb U_j\ee\bb U_{r_1}$.

\ 

$(\ga)$ \ $|A_1/E_L|<M_4$ \ or\  $|A_1/E_R|<M_4$.

\\Otherwise, suppose $\lan\bb X_1,\bb X_2,\ldots,\bb X_{M_4-1}\ran\sb
\{\bb U_i: i\in A_1\}$ is a sequence of pairwise $E_L$-nonequivalent 
representatives and that $\lan\bb Y_1,\bb Y_2,\ldots,\bb Y_{M_4-1}\ran
\sb\{\bb U_i: i\in A_1\}$ are pairwise $E_R$-nonequivalent. 
By $(\al)$ we know that for every $i,j<M_4$ there is some 
$h(i,j)<i^*$ with $\bb X_i\ee\bb Y_j\sim\bb U_{h(i,j)}$. Define a colouring 
$f\colon[M_4]^3\to\{0,1,\ldots,6\}$ by:
$$f(i,j,k)=\cases{0&if\/ $h(i,i)=h(j,k)$\cr 1&if\/ $h(i,i)=h(k,j)$\cr 
2&if\/ $h(j,j)=h(i,k)$\cr 3&if\/ $h(j,j)=h(k,i)$\cr 4&if\/ $h(k,k)=h(i,j)$\cr 
5&if\/ $h(k,k)=h(j,i)$\cr 6&otherwise.\cr}$$
(If more then one of these cases occurs, $f$ takes the minimal value.)

\\By the definition of\/ $M_4$ there is $B\sb M_4$ with $|B|=M_3$ \st $B$ is 
homogeneous with respect to $f$ and we let $f\res_B\equiv m$. Is it possible 
that $m<6$? Suppose for example that $m=0$, and choose $i<j<j_1<k$ from $B$.
If\/ $f(i,j,k)=0=f(i,j_1,k)$ we have $h(i,i)=h(j,k)$ and $h(i,i)=h(j_1,k)$. 
Hence $\bb X_i\ee\bb Y_i\sim\bb X_j\ee\bb Y_k$ and 	 
$\bb X_i\ee\bb Y_i\sim\bb X_{j_1}\ee\bb Y_k$. It follows that
$\bb X_j\ee\bb Y_k\sim\bb X_{j_1}\ee\bb Y_k$ and hence $\bb X_jE_L\bb X_{j_1}$
and this is impossible. 
The other five possibilities are eliminated similarly and we conclude that 
$$f\res_B\equiv 6.$$  
Let $A_2:=\{l<i^*: (\ex i\in B)(h(i,i)=l)\}$ and 
$A_3:=\{l<i^*: (\ex i\ne j\in B)(h(i,j)=l)\}$. By the choice of\/ $B$ and the 
above we have $A_2\cap A_3=\em$.
Note that $|A_2|\le|B|=M_3<K$ and $|A_3|\le|B|^2=(M_3)^2<K$. Hence by the 
$K$-randomness of\/ $\G$ there is some $k<i^*$ \st 
$$\big[l\in A_2 \ \imp \ C\mo R(\bb U_k,\bb U_l,\bb W)\big]\ \&\ 
\big[l\in A_3 \ \imp \ C\mo\neg R(\bb U_k,\bb U_l,\bb W)\big]$$
that is (as $R$ respects $\sim$)
$$(\dag) \ \ \ i\ne j\in B \ \ \ \imp \ \ \ C\mo [R(\bb U_k,\bb X_i\ee\bb Y_i,\bb W)\ \&\ 
\neg R(\bb U_k,\bb X_i\ee\bb Y_j,\bb W)].$$
By the definition of\/ $M_3=|B|$ we have $i\ne j\in B$ with 
$$(*) \ \ \ \th^n(R;\bb U_k,\bb Y_i,\bb W) = \th^n(R;\bb U_k,\bb Y_j,\bb W).$$
But $\th^n(C;\bb U_k,\bb X_i\ee\bb Y_i,\bb W) = 
\th^n(L;\bb U_k,\bb X_i,\bb W)+\th^n(R;\bb U_k,\bb Y_i,\bb W) =_{{\rm by\ }(*)}$

\\$\th^n(L;\bb U_k,\bb X_i,\bb W)+\th^n(R;\bb U_k,\bb Y_j,\bb W) = 
\th^n(C;\bb U_k,\bb X_i\ee\bb Y_j,\bb W)$.  Therefore:
$$C\mo R(\bb U_k,\bb X_i\ee\bb Y_i,\bb W) \iff 
C\mo R(\bb U_k,\bb X_i\ee\bb Y_j,\bb W)$$
and this is a contradiction to $(\dag)$, so $(\ga)$ is proved.

\ 

To conclude, assume $|A_1/E_L|<M_4$. Then, by 3.3 and as $\{u_i:i\in A_1\}$ is  
$(K,{K\over(M_2)^2})$-big, there is $A\sb A_1$ \st $\lan u_i: i\in A\ran$ 
is $(K+{K\over{(M_2)^2}},{K\over{(M_2)^2}\cdot M_4})$-big and \st for every 
$i,j\in A$, $\bbu_iE_L\bbu_j$. Fix $k^*\in A$, and define a sequence 
$\lan\bb V_i: i\in A\ran$ by: 
$$\bb V_i\res_L=\bb U_k^*\res_L \ {\rm and}\ \ \bb V_i\res_R=\bb U_i\res_R.$$
We want to show that for every $i\in A$ we have $\bb V_i\sim\bb U_i$.  
Indeed, as $\bbu_k^*E_L\bbu_i$ and by $(\be)$ we know that for every 
$r\in A_1$, $\bb U_{k^*}\ee\bb U_r\sim\bb U_i\ee\bb U_r$ and choosing $r=i$ we 
get $\bb U_{k^*}\ee\bb U_i\sim\bb U_i\ee\bb U_i$ i.e. $\bb V_i\sim\bb U_i$.
Hence $\lan v_i: i\in A\ran$ is $(K_1,K_2)$-big and all the $\bb V_i$'s 
coincide outside $R$. Hence $R$ is $(K_1,K_2)$-major.

By a similar argument we get: $|A_1/E_R|<M_4$ implies $L$ is $(K_1,K_2)$-major.  
\vv \qed

\

\pd Notation. $K_1$ and $K_2$ will be from now on the numbers from lemma 
3.9 above. 
\vv

\

The computations below will be useful the in following stages. For the moment 
assume that $\G$ is finite. 

Let $(L,R)$ be a Dedekind cut of\/ $C$ and $K>(M_3)^2$ as before. First note 
that as $\G$ is $K$-random we have
$$|\G|=\#(C)\ge 2^{2(K-1)}$$
By 3.9 we may assume that $L$ is $(K_1,K_2)$-major where 
$K_1=K+{K\over{(M_2)^2}}$ and $K_2={K\over{(M_2)^2} M_4}$. 
(The case $R$ is $(K_1,K_2)$-major is symmetric.) If\/ $K$ is big 
enough we get 
$$2^{2(K-1)}-K_1=2^{2(K-1)}-(K+{K\over{(M_2)^2}})>K>K_2 =
{K\over{(M_2)^2} M_4}$$
and by the definition of\/ $(K_1,K_2)$-major 
$$\#(L)\ge 2^{K_2} = 2^{{K\over{(M_2)^2} M_4}}.$$
By 3.7 $R$ is not $M_1(K_1+K_2)+1$-fat, i.e., 
$$\#(R)\le M_1(K_1+K_2) =
M_1(K+{K\over{(M_2)^2}}+{K\over{(M_2)^2} M_4})\le 2M_1 K$$
it follows that 
$$(*) \ \ \#(L)/(\#(R)+1)^2\ge 2^{K_2}/{(4(M_1)^2K^2+4(M_1)K+1)}=
2^{{K\over{(M_2)^2} M_4}}/{(4(M_1)^2K^2+4(M_1)K+1)}.$$

\

\pt Conclusion 3.10. For every $l<\om$ there is $K^*=K^*(l,n,d)<\om$ \st 
under the context in 3.1, if 

$K\ge K^*$,  

$(L,R)$ is a Dedekind cut of\/ $C$, 

$M=M(K,n,d)$ denotes the bouquet size of the major segment, 

$m=m(K,n,d)$ denotes the bouquet size of other segment, 

\\then $M/{(m+1)^2}>l\cdot K^2$.
\vv

\pd Proof. By the inequality $(*)$ above and noting that 
$M_1$, $M_2$, $M_3$ and $M_4$ do not depend on $K$.
\vv\qed

\

\pd Remark. By 3.7, if\/ $K$ is big enough then the segment that is not 
$(K_1,K_2)$-major is minor. We will always assume that.
\vvv

\

If we assume that the interpreted graph $\G$ is infinite then we can say 
that, if\/ $K$ is big enough, one segment will have an infinite bouquet size 
while the other will have an a priori bounded bouquet size. 

\pt Lemma 3.11. For every $n,d<\om$ there is $K^*=K^*(n,d)<\om$ \st if
$\II$, of dimension $d$ and depth $n$, is an interpretation of an 
infinite $K$-random graph $\G$ on $C$ and $K>K^*$, 

\\then there is $m<\om$, that depends only on $K$, $n$ and $d$, \st
if\/ $(L,R)$ is a Dedekind cut of\/ $C$:

$L$ (or $R$) has an infinite bouquet size,

$R$ (or $L$) has bouquet size that is at most $m$.
\vv

\pd Proof. By lemma 3.9 (letting $K^*=(M_3)^2$) we get that one of the 
segments is $(K_1,K_2)$-major and hence has infinite bouquet size. From 3.7 we get
that the other segment is not $K_3$-fat where $K_3$ depends only on $K$ and  
the interpretation $\II$ (i.e. $n$ and $d$). The required $m$ is that $K_3$,
which is good for every Dedekind cut.
\vv\qed
\sec\sec
\hh 4. Semi--Homogeneous subsets \par 

Our next step towards reaching a contradiction is of a combinatorial nature. 
In this section we introduce the notion of a \sg subset and show that the gap 
between the size of a set and the size of a \sg subset is reasonable.

\pd \dd 4.1.  Let $K,c<\om$, let $I$ be an ordered set and 
$f\colon{[I]^2}\to\{0,\ldots,c\}$. \par 
1) We call $T\sb I$ {\sl right \sg in I (for $f$ and $K$)} if for every 
$i<i^*$ from $T$ we have 
$|\{j\in I : j>i,\ f(i,j)=f(i,i^*)\}|\ge K$. \par
2) We call $T\sb I$ {\sl left \sg in I (for $f$ and $K$)} if for every 
$i<i^*$ from $T$ we have 
$|\{j\in I : j<i^*,\ f(j,i^*)=f(i,i^*)\}|\ge K$. \par
3) We call $T\sb I$ {\sl \sg in I (for $f$ and $K$)} if\/ $T$ is both right 
\sg and left semi-homogeneous. \par
4) We call $T\sb I$ {\sl right--nice [left--nice]} for $S\sb I$ 
(and for $f$ and $K$) if 

\\Max($S$)$<$Min($T$) [Min($S$)$>$Max($T$)] and for every $j\in T$, 
$S\cup\{j\}$ is right \sg [left semi-homogeneous] in $T\cup S$. 
\vvv

\

\pt Lemma 4.2. Let $K,c,I,f$ be as above.  Suppose $|I|>c\cdot N\cdot K$. 
Then, there is a right \sg subset $S\sb I$ of cardinality $N$. 
\vv 

\pd Proof. Let $i_0$ be Min($I$) and $T_0\sb I$ be of cardinality 
$\ge |I|-c\cdot(K-1)$ \st $T_0$ is right--nice for $\{i_0\}$,  
(just throw out every $j\in I$ \st $f(i_0,j)$ occurs less than $K$ times, 
there being at most $c\cdot(K-1)$ such $j$'s). 
Let $i_1$ be Min($T_0$) and $T_1\sb T_0$ be right--nice  
for $i_1$ of cardinality $\ge |I|-2c\cdot(K-1)$, (use the same argument).
Define $S_1:=\{i_0,i_1\}$. Clearly, for every $j\in T_1$, \ 
$S_1\cup\{j\}$  is right \sg in $I$. 

Proceed to define $i_2$ (=Min($T_1$) ), $T_2$, $S_2$ and so on. After 
defining $S_{N-2}$ and $T_{N-2}$ we have thrown out $(N-1)\cdot c\cdot(K-1)$ 
elements and as $|I|>c\cdot N\cdot K$ we can define $T_{N-1}, i_{N-1}$ and 
$S_{N-1}$ which is the required right \sg subset.  
\vv \qed

\pt Lemma 4.3. Let $K,c,I,f$ be as above.  
Suppose $|I|>c^2\cdot N\cdot K^2 \ =\ c\cdot (cNK)\cdot K$. \ 
Then, there is a \sg subset $T\sb I$ of cardinality $N$. 
\vv 

\pd Proof. Repeat the construction in the previous lemma to get $T^*\sb I$, 
right \sg in $I$ of cardinality $\ge c\cdot N\cdot K$ and now take 
$T\sb T^*$ left \sg in $T$ of cardinality $\ge N$. \par
\\$T$ is \sg in $I$.  
\vv \qed

\

We return now to the previous section, and its context. Recall that, given a 
cut $(L,R)$  we denoted by $M=M(K,n,d)$ the bouquet size of the major segment
and by $m=m(K,n,d)$ the bouquet size of the minor segment.

\pt Conclusion 4.4. In the context 3.1, for every $c,N<\om$ there is  
$K<\om$ \st if\/ $|I|\ge M(K,n,d)$ then for every 
$f\colon{[I]^2}\to\{0,\ldots,c\}$ there is a \sg subset of\/ $T\sb I$, 
for $f$ and $m(K,n,d)+1$, with $|T|\ge N$. 
\vv

\pd Proof. By lemma 4.3. we just need to ensure that 
$c^2\cdot N\cdot(m(K,n,d)+1)^2<M(K,n,d)$.  i.e.  
$M/{(m+1)^2}>c^2\cdot N$. This holds by conclusion 3.10.
\vv \qed
\sec\sec
\hh 5. The forcing \par   
The universe $V^P$ where no uniform interpretation exists is the same as in 
\ls. The forcing $P$ adds generic semi--clubs to each regular cardinal 
$>\ale$.

\pd Context. $V\models $ GCH 
\vvv

\pd \dd 5.1. Let $\la>\ale$  be a regular cardinal
\item{1)} $SC_\la$ := 
$\big\{ f: \ \ f\colon\al\to\{0,1\},\ \al<\la,\ \cf(\al)\le\ale\  \big\}$ \ 
where each $f$, considered to be a subset of\/ $\al$ (or $\la$), is a \sc. 
The order is inclusion. (So $SC_\la$ adds a generic \sc\ to $\la$). 
\item{2)} $Q_\la$ will be an iteration of the forcing $SC_\la$ with
length $\la^+$ and with support $<\la$. 
\item{3)} $P$ := $\lan P_\mu$, \qmu : $\mu {\rm \ a\ cardinal >\ale\ } 
\ran$ where \qmu \ is forced to be $Q_\mu$ if\/ $\mu$ is regular,
otherwise it is $\em$. The support of\/ $P$ is Easton's: each condition $p\in P$ 
is a function from the class of cardinals to names of conditions where 
the class $S$ of cardinals that are matched to non-trivial names is a set.
Moreover, when $\kappa$ is an inaccessible cardinal, $S\cap\kappa$ has 
cardinality $<\kappa$. 
\item{4)} $P_{<\la}, P_{>\la}, P_{\le\la}$ \ are defined naturally. 
For example $P_{<\la}$ is $\lan P_\mu$, \qmu : $\ale<\mu<\la \ran$. 
\vvv

\pd Discussion 5.2. Assuming GCH it is standard to see that 
$Q_{\la}$ satisfies the $\la^+$ chain condition and that $Q_{\la}$ and 
$P_{\ge\la}$ do not add subsets of\/ $\la$ with cardinality $<\la$. 
Hence, $P$ does not collapse cardinals and does not change cofinalities, so 
$V$ and $V^P$ have the same regular cardinals. 

Moreover, for a regular $\la>\ale$ we can split the forcing into 3 parts, 
$P=P_0*P_1*P_2$ where $P_0$ is $P_{<\la}$, $P_1$ is a $P_0$-name of the 
forcing $Q_\la$ and $P_2$ is a $P_0*P_1$-name of the forcing $P_{>\la}$ 
\st $V^P$ and $V^{P_0*P_1}$ have the same $H(\la^+)$.  

In the next sections, when we restrict ourselves to $H(\la^+)$ it will 
suffice to look only in $V^{P_0*P_1}$. 
\vv
\sec\sec
\hh 6. The contradiction (reduced case) \par

Collecting the results from the previous sections we will reach a 
contradiction from the assumption that (for a sufficiently large $K$), 
the monadic theory of some chain $C$ in $V^P$, interprets a radom graph 
$\G\in\Gamma_K$.

As we saw in section 3, an interpretation has a major segment. We will show  
below that there is a minimal one (and without loss of generality the 
segment is an initial segment). In this section we restrict ourselves to a 
special case: we assume that the minimal major initial segment is the whole 
chain $C$. Moreover, the chain $C$ is assumed to be regular cardinal $>\ale$. 

In the next section we will dispose of these special assumptions. However, 
the skeleton of those proofs will be the same as in this reduced case.

\

\pd \dd 6.1. Assume that $(C,<)$ interprets $\G\in\Ga_K$ by $\II$.
$D\sb C$ is a {\sl minimal $(K_1,K_2)$-major initial segment for $\II$} if 
$D$ is an initial segment of\/ $C$ which is a $(K_1,K_2)$-major segment and  
no proper initial segment $D'\sbb D$ is $(K_1,K_2)$-major. 
\vvv

\

\pt Fact 6.2. Suppose that $(C,<)$ interprets $\G\in\Ga_K$ by $\II$, where 
$K$ and $\II$ satisfy the assumption of lemma 3.9. Then there is 
a chain $(C^*,<^*)$ that interprets $\G$ by some $\II^*$ having the same 
dimension and depth as $\II$, \st there is $D^*\sb C^*$ which is  
a minimal $(K_1,K_2)$-major initial segment for $\II^*$. ($K_1$ and $K_2$ 
are as in lemma 3.9).
\vv
\pd Proof. (By \gu\ lemma 8.2). Let $L$ be the union of all the initial 
segments of\/ $C$ that are $(K_1,K_2)$-minor (note that if\/ $L$ is minor and 
$L'\sb L$ then $L'$ is minor as well). If L is $(K_1,K_2)$-major then set 
$D=L$, $C^*=C$, $\II^*=\II$ and we are done. 

Otherwise, let $D$ = $C\sm L$, by lemma 3.9 $D$ is major. Now if there is 
a proper final segment $D'\sbb D$ which is $(K_1,K_2)$-major then $C\sm D'$ 
is minor. But $(C\sm D')\supset L$, so that is impossible by maximality of 
$L$. Therefore $D$ is a minimal $(K_1,K_2)$-major (final) segment. Now take 
$C^*$ to be the inverse chain of\/ $C$. Clearly $D$ is a minimal 
$(K_1,K_2)$-major initial segment for an interpretation $\II^*$ of\/ $\G$ 
(that is obtained be replacing `$<$' by `$>$' in $\II$) having the same depth 
and dimension.  
\vv \qed

\

\\\pd Sketch of the proof.  Fixing an interpretation $\II$ (rather its depth and
dimension) we are trying to show that if\/ $K$ is large enough then in $V^P$ 
no chain $C$ interprets some $\G\in\Ga_K$ by $\II$. Towards a contradiction 
we choose $K$ \st  
$$\sqrt K>N_0>N_1>N_2>N_3>N_4>N_5>N_6$$ 
with: \par

\\(1) $N_6=\max\{2,n_1,n_2,n_3\}+1$ \ ($n_1,n_2,n_3$ are defined in 
assumption 5 below). 

\\(2) $N_5\to(N_6)^3_{32}$  \ i.e. a set of size $N_5$ has a homogeneous 
subset of size $N_6$ for colouring triplets into 
$32$ colours  (exists by Ramsey theorem).

\\(3) $N_4 = n_1\cdot N_5$.  

\\(4) $N_3=2\cdot N_4$.

\\(5) $N_2\to(N_3)^2_{n_3}$  \ (exists by Ramsey theorem).

\\(6) $N_1\to(N_2)^3_{32}$  \ (exists by Ramsey theorem).

\\(7) $N_0 = n_1\cdot N_1$. 

We start with a sequence $\lan \bb U_i : i<M \ran$\  of representatives 
for the elements of $\G$ ($M=\#(D)$ i.e. the bouquet size of a minimal 
major segment, in our case it is $\,\#(C)=|\G|$, possibly infinite), and 
gradually reduce their number until we get pairs that will satisfy 
(for a suitable \sc\ $a$ and a club $J$):
$$ i<j \ \imp \ [\bb U_i,\bb U_j]_a^J\sim \bb U_i .$$ 
These will be achieved at steps $1,2,3$. In steps $4,5$ we will get also:
$$ [\bb U_i,\bb U_j]_a^J\sim [\bb U_j,\bb U_i]_a^J .$$
Contradiction will be achieved when we show that some $[\bb U_i,\bb U_j]_a^J$ 
represents two different elements.
\vvv

\

\pd Assumptions. Our assumptions towards a contradiction are as follows: 

\\1. $(C,<)\in V^P$ interprets $\G\in\Ga_K$ by $\II=\lan U(\bb X,\bb W), 
E(\bb X,\bb Y,\bb W), R(\bb X,\bb Y,\bb W) \ran$, 

\\ \ \ \ $\lg(\bb X),\lg(\bb Y)$ and w.l..o.g\/ $\lg(\bb W)=d$, $n(\II)=n$. 

\\2. $C$ itself is the minimal $(K_1,K_2)$-major initial segment for $\II$.   
Moreover, $C=\la$, a regular cardinal $>\ale$. For every proper 
initial segment $D\sbb C$ we have $\#(D)<K_3$.  ($K_1$, $K_2$ and $K_3$ are 
from $\S3$, they depend only on $K$, $n$ and $d$). 

\\3. $m(*)=m(*)(n+d,4d)$ is as in the preservation theorem 1.11. 

\\4. $J=\lan\al_i:i<\la\ran\sb\la$ is an $m(*)$--suitable club for all the  
representatives that will be shuffled (there are only finitely many). $a\sb\la$ 
is a \sc, generic with respect to every relevant element including $J$ 
(again, finitely many), and see a remark later on.

\\5. $n_1$, $n_2$ and $n_3$ are defined as the number of possibilities for 
the following theories ($m(*)$ is as above): 

\centerline{$n_1 := |\bigl\{ \ath^{m(*)}(C;\bb X,\bb Y) : \ 
\bb X,\bb Y\sb C,\ \lg(\bb X),\lg(\bb Y)=d \bigr\}|$}

\centerline{$n_2 := |\bigl\{\ath^{m(*)}(C;\bb X,\bb Y,\bb Z) :
\ \bb X,\bb Y,\bb Z\sb C,\ \lg(\bb X),\lg(\bb Y),\lg(\bb Z)=d\bigr\}|$}

\centerline{$n_3 := | \bigl\{ \ath^{m(*)}(C;\bb X,\bb Y,\bb Z,\bb U) : \ 
\bb X,\bb Y,\bb Z,\bb U\sb C,\ 
\lg(\bb X),\lg(\bb Y),\lg(\bb Z),\lg(\bb U)=d \bigr\}|$}

\\6. \ $\sqrt K > N_0$. In addition, $K$ is large for 
$l:=(n_2)^2\cdot N_0\cdot(2|T_{n,3d}|+1)^2$ as in conclusion 3.10 i.e. 
$M/{(m+1)^2}>l\cdot K^2$ (this is possible as $l$ depends only on $n(\II)$ 
and $d(\II)$).
\vvv 

\

\\To get started we need another observation that does not depend on the 
special assumption on the minimal major segment. 

\pd \dd 6.3. Suppose $D$ is the minimal $(K_1,K_2)$-major initial segment for 
the interpretation. The {\sl vicinity} of a representative $\bb X$ denoted by 
$[\bb X]$ is the collection of representatives 

\\$\{ \bb Y :\ $some \ $\bb Z\sim\bb Y$\ coincides with\ $\bb X$\  
outside\ some proper (hence minor) initial segment of\ $D \}$. 
\vvv

\

\pt Lemma 6.4.  (1) Every vicinity $[\bb X]$ is the union of at most 
$m=m(K,n,d)$ \ (the bouquet size of a minor segment) different equivalence 
classes.   \par
\\(2) From $\th^{n+d}(D;\bb U_1,\bb U_2,\bb W)$ we can compute the truth 
value of:  ``$\,\bb U_1$ is in the vicinity of\/ $\bb U_2$''. 
\vvv

\pd Proof. If\/ $(1)$ does not hold then there is a proper initial  
segment $D'$ of the minimal major initial segment $D$ with \/ 
$\#(D')>m$ which is impossible. $(2)$ is clear.
\vv \qed

\

We are ready now for a contradiction: 

{\bf STEP 1:} \ Let $N_0\ldl N_6$ be as above and let $K$ be as in assumption 
6.  

\\Let $\lan \bb U_i : i<M \ran$, be a list of representatives for the 
elements of\/ $\G\in\Ga_K$ that is interpreted by $\II$ on $C$.
Let $f$ be a colouring of\/ $[M]^2$ into $n_2$ colours defined by  
$$f(i,j):=\ath^{m(*)}(C;\bb U_i,\bb U_j,\bb W).$$ 

We would like to get a \sg subset of\/ $M$ for $f$ and $m+1$ of size $N_0$. 
If\/ $\G$ is finite then this is possible by assumption 6 and conclusion 4.4.
Of course if\/ $\G$ is infinite (i.e. $M\ge\ale$) we can even get a 
homogeneous one. 

Let then $S'\sb\{0\ldl M-1\}$ be \sg and look at $B':=\lan\bbu_i:i\in S'\ran$.
As $N_0=n_1\cdot N_1$ we can choose 
$$B:=\lan\bbu_i:i\in S\ran$$ 
\st $S\sb S'$ is of size $|N_1|$ and \st $\ath^{m(*)}(C;\bbu_i,\bbw)$ is 
constant for every $i\in S$.

\

{\bf STEP 2:} \ We start shuffling the members of\/ $B$ along $a$ and $J$. 
Note that by the choice of\/ $B$ and $m(*)$ and by the preservation theorem 
$$i,j\in S \imp \th^n(C;\bbu_i,\bbw)=\th^n(C;[\bbu_i,\bbu_j]_a^J,\bbw).$$
It follows that the results of the shufflings are representatives as well, 
that is
$$i,j\in S \imp C\mo U([\bb U_i,\bb U_j]_a^J,\bbw).$$
Define for $i<j\in S$
$$k(i,j) := \min\Big\{ k : \big(k\in S\ \&\ [\bb U_i,\bb
U_j]_a^J\sim\bb U_k\big)\vee\big(k=M\big)\Big\}$$
By the choice of\/ $N_1$ there is a subset $A\sb S$, of size $N_2$, \st for 
every $\bb U_i,\bb U_j,\bb U_l$ with $i<j<l$ and $i,j,l\in A$, the following 
five statements have a constant truth value:

$k(j,l)=i$,

$k(i,l)=j$,

$k(i,j)=i$,

$k(i,j)=j$,

$k(i,j)=l$.

\\Moreover, if there is a pair $i<j$ in $A$ \st $k(i,j)\in A$ then:

$${\rm either\ for\ every\ } i<j\ {\rm from\ } A, \ k(i,j)=i\ {\rm or\ for\ 
every\ } i<j\ {\rm from\ } A, \ k(i,j)=j.$$

\\The reason is the following: suppose that $k(\al,\be)=\ga\in A$ for some 
$\al<\be$ from  $A$. If\/ $\ga<\al$ then $k(j,l)=i$ for all $i<j<l$ from $A$
but $k$ is one valued. Similarly, the possibilities $\al<\ga<\be$ and 
$\be<\ga$ are ruled out. We are left with $\ga=\al$ or $\ga=\be$ and 
apply homogeneity. 

\

{\bf STEP 3:} \ The aim now is to find a pair $i<j$ from $A$ with 
$k(i,j)\in A$. Define $A^*$ to be the results of the shufflings:
$$A^*:=\big\{k: (\ex i<j\in A)\big([\bb U_i,\bb U_j]_a^J\sim\bb U_k\big)\big\}$$
and it is enough to show that $A^*\cap A\not=\em$.

If not, as $|A|=N_2<K$ and $|A^*|\le|A|^2<K$ (we chose $\sqrt K > N_0$),
there is a representative $\bb V_A$ \st 
$$\bigwedge_{i\in A}[C\mo R(\bbu_i,\bbv_A,\bbw)]\wedge
\bigwedge_{i\in A^*\sm A}[C\mo \neg R(\bbu_i,\bbv_A,\bbw)].$$
As $N_2>n_2$ there is $i<j \in A$ with:
$$\ath^{m(*)}(C;\bb U_i,\bb V_A,\bb W)=\ath^{m(*)}(C;\bb U_j,\bb V_A,\bb W)$$
and by the preservation theorem 
$$ \ \ (*) \ \
\th^n(C;[\bb U_i,\bb U_j]_a^J,\bb V_A,\bb W)=\th^n(C;\bb U_i,\bb V_A,\bb W).$$
Now, $[\bb U_i,\bb U_j]_a^J\sim\bb U_k$ for some $k\in A^*$ but by $(*)$
$$C\mo R([\bbu_k,\bb V_A,\bb W).$$
Therefore, by the choice of\/ $\bb V_A$, we have $k\in A$. It follows that 
$k\in A\cap A^*$ so $A^*\cap A\not=\em$ after all. 

\\The aim is fulfilled and we may assume w.l.o.g that 
$$i<j\in A \imp [\bb U_i,\bb U_j]_a^J\sim\bb U_i.$$

\

{\bf STEP 4:} \ The aim now is to show that 
$$\otimes \ \ i<j \in A \  \imp [\bb U_i,\bb U_j]_a^J\sim [\bb U_j,\bb U_i]_a^J \
\ \big(= [\bb U_i,\bb U_j]^J_{\la\sm a}\big).$$
Returning to the discussion in $\S5$, we have mentioned so far only a finite 
number of elements from $H(\la^+)^{V^P}$, (including $J$). Everything  
already belongs to $H(\la^+)^{V^{P_0*P_1}}$ where $P_0$ is $P_{<\la}$ and 
$P_1$ \/ is a $P_0$-name  for $Q_\la$ which is an iteration of length 
$\la^+$ with support $<\la$  (we assume that the ground universe $V$ 
satisfies GCH). Moreover, an initial segment of\/ $P_0*P_1$, denoted by 
$P_0*P_1\res_\be$ adds all the relevant elements and we can choose the \sc\  
$a$ as the one that is generated in the $\be$'th stage of\/ $P_1$.

\\Let  $p\in P_0*P_1$ be a condition that forces all the statements about 
the representatives we mentioned so far (e.g. 
$i<j \in A \imp [\bb U_i,\bb U_j]_a^J\sim\bb U_i$).
We think about $p$ as a function with domain $\{-1\}\cup\la^+$ \st 
$p(-1)\in P_0$ and for $\al\in \la^+$, $p(\al)\in SC_\la$. under this 
notation $p(\be)$ is an initial segment of\/ $a$ and w.l.o.g a member of 
$V^{P_0}$ (and not a name for one). Let $\ga^*={\rm Dom}(p(\be))$. We may 
assume that $\cf(\ga^*)=\ale$. Let $\ga:=\al_{\ga^*}\in J$ (so $\cf(\ga)=\ale$ 
as well). 

By homogeneity of the forcing,  
$b:=(a\cap\ga)\cup\big[(\la\sm a)\cap[\ga,\la)\big]$ is a \sc\  of\/ $\la$ that 
is also generic with respect to the relevant elements. 
We denote from now on, for $\bbu,\bbv\sb\la$, 
$$\bb U\ee\bb V :=
\big(\bb U\cap\ga\big)\cup\big(\bb V\cap[\ga,\la)\big).$$
For proving $\otimes$ we will show that: 

$(\al)$ \ $[\bb U_i,\bb U_j]^J_b\sim\bb U_i$ for all $i<j$ from $A$,

$(\be)$ \ $[\bb U_i,\bb U_j]^J_{\la\sm a}\ee\bb U_k \ \sim \ \bb U_k$ for 
all $i,j,k$ from $A$,

$(\ga)$ \ $[\bb U_i,\bb U_j]^J_{\la\sm a}\ee\bb U_i \ \sim \
[\bb U_i,\bb U_j]^J_{\la\sm a}$ for all $i<j$ from $A$.

\

{\bf STEP 5:} \ Let's prove the claims: 

\\$(\al)$: By homogeneity of the forcing everything that $p$ forces for $a$ 
it forces for $b$.

\

\\$(\be)$: Recall that for every $i,j,k\in A$ we have
$$\ath^{m(*)}(C;\bb U_i,\bb W)=\ath^{m(*)}(C;\bb U_j,\bb W)=
\ath^{m(*)}(C;\bb U_k,\bb W).$$
As $m(*)=m(*)(n+d,4d)$ and $\ga\in J$ satisfies $\cf(\ga)=\ale$ 
we have by the second part of preservation theorem 1.11 
$$\th^{n+d}(C;[\bb U_i,\bb U_j]_a^J,\bb W)\res_{[0,\ga)}=
\th^{n+d}(C;\bb U_i,\bb W)\res_{[0,\ga)}=
\th^{n+d}(C;\bb U_j,\bb W)\res_{[0,\ga)}.$$
Similarly 
$$\th^{n+d}(C;\bb U_i,\bb W)\res_{[0,\ga)}=
\th^{n+d}(C;\bb U_k,\bb W)\res_{[0,\ga)}$$
and it follows that for every $i,j,k\in A$:
$$(\dagger) \ \ \th^{n+d}(C;[\bb U_i,\bb U_j]_a^J,\bb W)\res_{[0,\ga)}=
\th^{n+d}(C;\bb U_k,\bb W)\res_{[0,\ga)}.$$
Now by the composition theorem 
$$\th^{n+d}(C;[\bb U_i,\bb U_j]_a^J\ee\bb U_k,\bb W) =
\th^{n+d}(C;[\bb U_i,\bb U_j]_a^J,\bb W)\res_{[0,\ga)} + \th^{n+d}(C;\bb
U_k,\bb W)\res_{[\ga,\la)}$$
and this equals by $(\dagger)$ 
$$\th^{n+d}(C;\bb U_k,\bb W)\res_{[0,\ga)} +
\th^{n+d}(C;\bb U_k,\bb W)\res_{[\ga,\la)} =
\th^{n+d}(C;\bb U_k,\bb W).$$
As the theories are equal and as $\bb U_k$ is a representative, there is some 
$l<M$ (not necessarily in $A$) \st 
$$[\bb U_i,\bb U_j]_a^J\ee\bb U_k\sim \bbu_l.$$
If\/ $l=k$ everything is fine. Otherwise assume $l>k$ (symmetrically for $l<k$)
for a contradiction.

By the definition of vicinity we see that $\bb U_l\in[\bb U_k]$ 
and this is reflected in $\th^{n+d}(C;\bb U_k,\bb U_l,\bb W)$.
Now $k\in A\sb S$ and $S$ was chosen to be \sg in $M$. Therefore there are 
$l_0<l_1<\ldots<l_m<M$ with 
$$\bigwedge_{i<(m+1)} \th^{n+d}(C;\bb U_k,\bb U_{l_i},\bb W)=
\th^{n+d}(C;\bb U_k,\bb U_l,\bb W).$$
Hence
$$\bigwedge_{i<(m+1)}\big(U_{l_i}\in[U_k]\big)$$
but by 6.4 a vicinity contains at most $m$ pairwise nonequivalent 
representatives, a contradiction. We conclude that $l=k$.

Therefore, for every $i,j,k$ from $A$ we have
\/ $[\bb U_i,\bb U_j]_a^J\ee\bb U_k\sim \bbu_k$.
Substituting $i$ and $j$ we get: for every $i,j,k$ from $A$, \/ 
$[\bb U_j,\bb U_i]_a^J\ee\bb U_k \sim\bb U_k$ or:
$$[\bb U_i,\bb U_j]_{\la\sm a}\ee\bb U_k \ \sim \ \bb U_k.$$
This is claim $(\be)$.

\

\\$(\ga)$: Now suppose $i<j$ are from $A$. By definition, for every 
$\bb P\sb C$ the theory \/ $\lath^{m(*)}(C;\bbp)$ determines 
(and is determined by) $\ath^{m(*)}(C;\bbp)$. Therefore, 
$$\ath^{m(*)}(C;\bb U_i,\bb W)=\ath^{m(*)}(C;\bb U_j,\bb W) \ \& \
\lath^{m(*)}(C;\bb U_i,\bb W)=\lath^{m(*)}(C;\bb U_j,\bb W).$$
Applying the preservation theorem for $a$ and $\la\sm a$ we get 
$$\th^{n+d}(C;[\bb U_i,\bb U_j]_a^J,\bb W)\res_{[0,\ga)}=
\th^{n+d}(C;\bb U_i,\bb W)\res_{[0,\ga)}=
\th^{n+d}(C;\bb U_j,\bb W)\res_{[0,\ga)}$$
and
$$\th^{n+d}(C;[\bb U_i,\bb U_j]_{\la\sm a}^J,\bb W)\res_{[0,\ga)}=
\th^{n+d}(C;\bb U_i,\bb W)\res_{[0,\ga)}=
\th^{n+d}(C;\bb U_j,\bb W)\res_{[0,\ga)}$$
so
$$\th^{n+d}(C;[\bb U_i,\bb U_j]_a^J,\bb W)\res_{[0,\ga)}=
\th^{n+d}(C;[\bb U_i,\bb U_j]_{\la\sm a}^J,\bb W)\res_{[0,\ga)}.$$
Therefore, as $\th^{n+d}(C;\bbp)$ determines $\th^{n+d}(C;\bbp,\bbp)$:
$$(\ddagger) \
\th^{n+d}(C;[\bb U_i,\bb U_j]_a^J,[\bb U_i,\bb U_j]_a^J,\bb W)\res_{[0,\ga)}=
\th^{n+d}(C;[\bb U_i,\bb U_j]_{\la\sm a}^J,[\bb U_i,\bb U_j]_{\la\sm a}^J,
\bb W)\res_{[0,\ga)}.$$
By  $(\al)$ and $(\be)$ we know that
$$[\bb U_i,\bb U_j]_b\ \sim \ \bb U_i\ \sim \
[\bb U_i,\bb U_j]^J_{\la\sm a}\ee\bb U_i$$
and the equivalence is reflected by 
$\th^n(C;[\bb U_i,\bb U_j]^J_{\la\sm a}\ee\bb U_i,[\bb U_i,\bb U_j]_b,\bbw)$.
Clearly $\th^n$ is determined by $\th^{n+d}$. Hence:

$\th^n(C;[\bb U_i,\bb U_j]_a^J\ee\bb U_i,[\bb U_i,\bb U_j]_b^J)=$
\ (by $(\ddagger)$)

$\th^n(C;[\bb U_i,\bb U_j]_a^J,[\bb U_i,\bb U_j]_a^J,\bb W)\res_{[0,\ga)}+
\th^n(C;\bb U_i,[\bb U_i,\bb U_j]_{\la\sm a}^J,\bb W)\res_{[\ga,\la)}=$

$\th^n(C;[\bb U_i,\bb U_j]_{\la\sm a},[\bb U_i,\bb U_j]_{\la\sm a},\bb W)
\res_{[0,\ga)}+
\th^n(C;\bb U_i,[\bb U_i,\bb U_j]_{\la\sm a},\bb W)\res_{[\ga,\la)}=$

$\th^n(C;[\bb U_i,\bb U_j]_{\la\sm a}\ee\bb U_i,[\bb U_i,\bb
U_j]_{\la\sm a}).$

\\Therefore 
$$[\bb U_i,\bb U_j]^J_{\la\sm a}\ee\bb U_i \ \sim \
[\bb U_i,\bb U_j]_{\la\sm a}^J$$
and $(\ga)$ is proved.

\

>From $(\be)$ and $(\ga)$ we conclude 
$$[\bb U_i,\bb U_j]_{\la\sm a}^J \ \sim \ \bb U_i \ \sim \
[\bb U_i,\bb U_j]_a^J$$
and $\otimes$ is proved.

\

{\bf STEP 6:} \ By definition of\/ $N_2=|A|$, $N_3$ and $N_4$ there is a  
sub-sequence of\/ $\lan \bb U_i : i\in A \ran$ that will be denoted for 
convenience (while preserving the order between the indices) by 
$\lan \bbh_i : i<N_3=2N_4 \ran$ \st for every $i<j<2N_4$ and $r<l<2N_4$:

$(i)$ \/ $\ath^{m(*)}(C;\bbh_i,\bbh_j,\bb W) =
    \ath^{m(*)}(C;\bbh_r,\bbh_l,\bb W)$ (by defining a colouring of pairs 
from $N_2$).

$(ii)$ \/ $[\bbh_i,\bbh_j]_a^J\sim[\bbh_i,\bbh_j]_{\la\sm a}\sim\bbh_i$
(by steps 3 and 5).

\\For $i<N_4$ let $\bbq_i$ a representative that satisfies 
$$\bigwedge_{\al\in[i,2N_4-i)}\Big(C\mo R(\bbh_\al,\bbq_i,\bbw)\Big)\wedge
\bigwedge_{\al\in[0,i)\cup[2N_4-i,2N_4)}\Big(C\mo\neg
R(\bbh_\al,\bbq_i,\bbw)\Big).$$
As $N_4$ is big enough there is $T\sb\{0\ldl N_4-1\}$ with $|T|=N_6$ \st
if\/ $i,j\in T$ then either $[\bbq_i,\bbq_j]_a^J\sim\bbq_i$ or 
$[\bbq_i,\bbq_j]_a^J\sim\bbq_j$.
To get $T$ repeat steps 1, 2 and 3 while substituting 
$\lan\bbu_i:i<S'\ran$ by $\lan\bbq_i:i<N_4\ran$, and $N_0$, $N_1$, $N_2$ by
$N_4$, $N_5$ and $N_6$ respectively.
Note that we lose generality by chosing one of the possibilities. 

Now choose $i,j\in T$ (by  $N_6>n_3$) \st
$$\ath^{m(*)}(C;\bbh_i,\bbh_{2N_4-i},\bbq_i,\bb W) =
\ath^{m(*)}(C;\bbh_j,\bbh_{2N_4-j},\bbq_j,\bb W)$$
and shuffle along $a$ and $J$:

$\th^n(C;\bbh_i,\bbh_{2N_4-i},\bbq_i,\bb W)=$

$\th^n(C;[\bbh_i,\bbh_j]_a^J,[\bbh_{2N_4-i},\bbh_{2N_4-j}]_a^J,
[\bbq_i,\bbq_j]_a^J,\bb W)=$

$\th^n(C;[\bbh_i,\bbh_j]_a^J,[\bbh_{2N_4-j},\bbh_{2N_4-i}]^J_{\la\sm a},
[\bbq_i,\bbq_j]_a^J,\bb W)$

\\but $[\bbh_i,\bbh_j]_a^J\sim\bbh_i$, and by step 5,
$$[\bbh_{2N_4-j},\bbh_{2N_4-i}]^J_{\la\sm a}\sim\bbh_{2N_4-j}.$$
Now from ``$[\bbq_i,\bbq_j]_a^J\sim \ \bbq_i$ or 
$[\bbq_i,\bbq_j]_a^J\sim \ \bbq_j$'' and the equality of the theories 
$\th^n$:
$$C\mo\big(R(\bbh_i,\bbq_i,\bbw)\ \&\ R(\bbh_{2N_4-j},\bbq_i)\big)$$
or
$$C\mo\big(\neg R(\bbh_i,\bbq_j,\bbw) \ \&\ \neg
R(\bbh_{2N_4-j},\bbq_j).$$
Both possibilities contradict the choice of the $\bbq_i$'s !

\

\pd First Remark. So $J$ and $a$ are chosen as follows:  getting 
$\lan\bbu_i:i\in S\ran$ at step 1 ($|S|=N_1$) choose for every subset $A\sb S$ 
a representative $\bbv_A$ that separates $\lan\bbu_i:i\in A\ran$ from 
$\lan\bbu_i:i\in S\sm A\ran$ (some of these will be the $\bbq_i$'s from step 6).
$J$ is an $m(*)$-suitable for all these elements and $a$ is a \sc\  that is 
generic with respect to all of these. Clearly, only finitely many 
elements are involved.

\

\pd Second Remark. Note that genericity was used only at stages 4 and 5 
(i.e. to prove $[\bbu_i,\bbu_j]^J_a\sim[\bbu_j,\bbu_i]^J_a$).

\

We proved the following:

\pt Theorem 6.5.  Let 
$\lan U(\bbx,\bbz),E(\bbx,\bby,\bbz),R(\bbx,\bby,\bbz)\ran$ be a sequence of
formulas of dimension $d$ and depth $n$. 

\\Then there is $K<\om$, that depends only on $d$ and $n$ such that, in $V^P$, 
for no chain $C$ and parameters $\bbw\sb C$:

$(i)$ \/ $C$ is isomorphic to a regular cardinal $\la>\ale$,

$(ii)$ \/ $\II=\lan U(\bbx,\bbw),E(\bbx,\bby,\bbw),R(\bbx,\bby,\bbw)\ran$
is an interpretation for some $\G\in\Ga_k$ in $C$,

$(iii)$ \/ $C$ is the minimal $(K_1,K_2)$-major initial (or final) segment 
for $\II$.
\vv\qed
\sec\sec
\hh 7. Generality                 \par 
Our aim in this section is to achieve full generality of the interpreting 
chain $C$ and its minimal initial major segment $D$. There are three stages: 

$(I)$ \ $D\sb C$, $D\ne C$ but $D$ is (isomorphic to) a regular 
cardinal $\la>\ale$. 

$(II)$ \ $D=C$, $D$ general. 

$(III)$ \ $C$ and $D$ are general.

\\Let us just remark that always $\cf(D)>\ale$, otherwise we can 
prove the non existence of interpretations even from ZFC. 

We will elaborate on stages (I) and (II), stage (III) is a simple 
combination of the techniques.

\

\pd Chopping off the final segment. We are trying now to get a contradiction 
from the same assumptions as in the previous section except for the following: 
the minimal $(K_1,K_2)$-major initial segment $D$ that is a regular cardinal 
is not necessarily equal to the interpreting chain $C$. $a$ and $J$ are 
therefore subsets of\/ $D$.

The basic idea of the proof is that if\/ $t^*$ is fixed and known in advance 
then to know $t_i+t^*$ all we need to know is $t_i$. Here $t_i$ are the 
restrictions of the information (partial theories) to $D$ and $t^*$ is 
the restriction to $C\sm D$ which can be assumd to be fixed, as many 
representatives coincide outside $D$. 

We do not specify the exact size of\/ $K$ (which should be slightly 
bigger than in the previous case). It should be apparent however that 
$K$ depends on $n$ and $d$ only and is obtained by repeated applications of 
the Ramsey functions. 

\

{\bf Preliminary Step:} Let $\lan\bbu_i:i<|\G|\ran$ a list of representatives 
for the elements of\/ $\G$. By definition of\/ $D$, we may assume that 
$\lan\bbu_i:i<M=\#(D)\ran$ is a list of representatives for a $(K_1,K_2)$-major 
subset of\/ $\G$ and all of them coincide outside $D$. 
Denote $D^@:=C\sm D$ and for $i<M$:

$\vvu_i:=\bbu_i\cap D$,

$\ubb:=\bbu_i\cap D^@$,

$\vvw:=\bbw\cap D$,

$\wbb:=\bbw\cap D^@$.

\

\pd \dd 7.1. 

(1) Define on $\P(D)$ a unary relation $U^*(\bbx)$ and binary relations 
$\bbx\sim^*\bby$ and $R^*(\bbx,\bby)$, with arity $d$ by: 

\centerline{$U^*(\bbx) \iff C\mo U(\bbx\cup\ubb,\bbw)$,}

\centerline{$\bbx\sim^*\bby \iff C\mo E(\bbx\cup\ubb,\bby\cup\ubb,\bbw)$,}

\centerline{$R^*(\bbx,\bby) \iff C\mo R(\bbx\cup\ubb,\bby\cup\ubb,\bbw)$.}

\\(When $i,j<M$ for instance, then $R^*(\vvu_i,\vvu_j)$ holds if and only if 
$C\mo R(\bbu_i,\bbu_j,\bbw)$).

(2) If\/ $i,j<|\G|$ and $\bbu_i\cap D^@=\bbu_j\cap D^@$ we denote  
$$[\bbu_i,\bbu_j]_a^J :=  
[\bbu_i\cap D,\bbu_j\cap D]_a^J\cup\big(\bbu_i\cap D^@\big)$$
\\(If\/ $i,j<M$ for instance then $[\bbu_i,\bbu_j]_a^J$ is 
$[\vvu_i,\vvu_j]_a^J\cup\ubb$).
\vvv

\

\pt Fact 7.2. $H(\la^+)^{V^P}$ computes correctly $\,\sim^*$, \/ 
$U^*$ \/ and \/ $R^*$ from $\ath^{m(*)}$
\vv
\pd Proof. Take for example $\sim^*$: $\bbx\sim^*\bby$ is determined by 
$\th^n(C;\bbx\cup\ubb,\bby\cup\bbu^@,\bbw)=$
$$\th^n(D;\bbx,\bby,\bbw^*)+\th^n(D^@;\ubb,\ubb,\wbb).$$
The second theory is fixed for every $\bbx,\bby\sb D$. Hence 
(e.g. by the finite number of possibilities) all we need to know is the first 
theory, which is computed correctly in $H(\la^+)^{V^P}$ from $\ath^{m(*)}$. 
\vv\qed

\

We proceed by immitating the previous proof:

{\bf STEP 1:} \ Define $B':=\lan \bb U_i:i\in S'\ran$ where $S'\sb\{0\ldl M-1\}$ 
is \sg, and 
$$B:=\lan\bb U_i:i\in S\ran$$
\st $S\sb S'$, $|S|$ finite and big enough, with   
$\ath^{m(*)}(D;\vvu_i,\vvw)$ constant for $i\in S$. 

\

{\bf STEP 2:} \ Shuffle the members of\/ $B$ along $a$ and $J$ as in definition 
7.1. Note that by the choice of\/ $B$ and the preservation theorem 
$$i,j\in S \imp \th^n(C;\bbu_i,\bbw)=\th^n(C;[\bbu_i,\bbu_j]_a^J,\bbw)$$
and therefore the resuts are representatives as well i.e. 
$$i,j\in S \imp C\mo U([\bb U_i,\bb U_j]_a^J,\bbw).$$
Define for $i<j\in S$
$$k(i,j) := \min\Big\{ k : \big(k\in S\ \&\ [\bb U_i,\bb
U_j]_a^J\sim\bb U_k\big)\vee\big(k=M\big)\Big\}$$
equivalently
$$k(i,j) := \min\Big\{ k : \big(k\in S\ \&\ [\vvu_i,\vvu_j
]_a^J\sim^*\vvu_k\big)\vee\big(k=M\big)\Big\}$$
Let $A\sb S$ be large enough, homogeneous with the colouring into 32 colours 
we used before.

\

{\bf STEP 3:} \ The aim is to find $i<j$ from $A$ with $k(i,j)\in A$. Let
$$A^*:=\big\{k<|\G|:(\ex i<j\in A)\big([\bb U_i,\bb U_j]_a^J\sim\bb U_k\big)\big\}$$
and let's show that $A^*\cap A\not=\em$.

Othewise, there is some $\bb V_A$ (not necessarilly from $\lan\bbu_i:i<M\ran$) 
that separates these two disjoint collections of representatives:
$$\bigwedge_{i\in A}[C\mo R(\bbu_i,\bbv_A,\bbw)]\wedge
\bigwedge_{i\in A^*\sm A}[C\mo \neg R(\bbu_i,\bbv_A,\bbw)].$$
We may assume that there are $i<j$ from $A$ with 
$$\ath^{m(*)}(C;\bb U_i,\bb V_A,\bb W)\res_D=
\ath^{m(*)}(C;\bb U_j,\bb V_A,\bb W)\res_D.$$
By the preservation theorem 
$$ \ \ (*) \ \
\th^n(C;[\bb U_i,\bb U_j]_a^J,\bb V_A,\bb W)\res_{D}=
\th^n(C;\bb U_i,\bb V_A,\bb W)\res_{D}$$
and in addition 
$$ \ \ (**) \ \
\th^n(C;[\bb U_i,\bb U_j]_a^J,\bb V_A,\bb W)\res_{D^@}=
\th^n(D^@;\ubb,\bbv_A\cap D^@,\wbb)=
\th^n(C;\bb U_i,\bb V_A,\bb W)\res_{D^@}.$$
Now $[\bb U_i,\bb U_j]_a^J\sim\bb U_k$ for some $k\in A^*$ but by $(*)$
and $(**)$ and the composition theorem:
$$C\mo R([\bb U_i,\bb U_j]_a^J,\bb V_A,\bb W).$$
Therefore, by the choice of\/ $\bb V_A$, $k\in A$. As $k\in A^*$ it follows 
that $A^*\cap A\not=\em$ after all.

As before we may conclude that, without loss of generality: 
$$i<j\in A \imp [\bb U_i,\bb U_j]_a^J\sim\bb U_i$$
equivalently (and this is known even by $H(\la^+)^{V^P}$)
$$i<j\in A \imp [\vvu_i,\vvu_j]_a^J\sim^*\vvu_i$$

\

{\bf STEPS 4,5:} \  We work inside $H(\la^+)^{V^P}$ and concentrate on 
$\lan\vvu_i:i\in A\ran$. The aim is to show that 
$$\otimes \ \ i<j \in A \imp [\vvu_i,\vvu_j]_a^J\sim^*[\vvu_j,\vvu_i]_a^J.$$
Let $p\in P_0*P_1$ be a condition that forces all the facts we showed so far 
about $\sim^*$, $U^*$ and $R^*$ and the $\vvu_i$'s such as 
$i<j \in A \imp [\vvu_i,\vvu_j]_a^J\sim\vvu_i$. 
As before we define a generic \sc\ $b\sb\la$ and show that:

$(\al)$ \  $[\vvu_i,\vvu_j]^J_b\sim^*\vvu_i$ for every $i<j$ from $A$,

$(\be)$ \  $[\vvu_i,\vvu_j]^J_{\la\sm a}\ee\vvu_k \ \sim^* \ \vvu_k$ for 
every $i,j,k\in A$,

$(\ga)$ \  $[\vvu_i,\vvu_j]^J_{\la\sm a}\ee\vvu_i \ \sim^* \
[\vvu_i,\vvu_j]^J_{\la\sm a}$ for every $i<j$ from $A$.

\\The proofs are exactly the same as in the previous section
(substituting $C$, $\bbw$ and $\sim$ by $D$, $\vvw$ and $\sim^*$), and
from these facts we can deduce $\otimes$. Leaving $H(\la^+)^{V^P}$ we find 
that what we proved in $V^P$ is:
$$\odot \ \ i<j \in A \imp [\bbu_i,\bbu_j]_a^J\sim[\bbu_j,\bbu_i]_a^J$$

\

{\bf STEP 6:} \  As $|A|$ is big enough we have a sub-sequence of 
$\lan \bb U_i:i\in A \ran$ that will be denoted for convenience 
(while preserving the order between the indices) by 
$\lan \bbh_i:i<N^*_3=2N^*_4\ran$ \st for every $i<j<2N^*_4$ and 
$r<l<2N^*_4$ ($N^*_4$ is a sufficiently big number as usual):

$(i)$ \/ $\ath^{m(*)}(D;\vvp_i,\vvp_j,\vvw) =
    \ath^{m(*)}(D;\vvp_r,\vvp_l,\vvw)$ 

$(ii)$ \/ $\th^n(D^@;\pbb_i,\pbb_j,\wbb) =
    \th^n(D^@;\pbb_r,\pbb_l,\wbb)$ \ ($\bbh_i\cap D^@$ is constant),

$(iii)$ \/ $[\bbh_i,\bbh_j]_a^J\sim[\bbh_i,\bbh_j]_{\la\sm a}\sim\bbh_i$.

\\For $i<N^*_4$ let $\bbq_i$ a representative (not necessarily from 
$\lan\bbu_i:i<M\ran$) that satisfies:
$$\bigwedge_{\al\in[i,2K_4-i)}\Big(C\mo R(\bbh_\al,\bbq_i,\bbw)\Big)\wedge
\bigwedge_{\al\in[0,i)\cup[{2N^*_4}-i,{2N^*_4})}\Big(C\mo\neg
R(\bbh_\al,\bbq_i,\bbw)\Big)$$
>From the $\bbq_i$'s extract $\lan\bbq_i:i\in T^1\ran$, with 
$|T^1|=N^*_5$ large enough \st for every $i<j$ from $T^1$:
$$\Big(\ath^{m(*)}(D,\vvq_i,\vvw)=\ath^{m(*)}(D,\vvq_j,\vvw)\Big) \ \&
\ \Big(\th^n(D^@,\qbb_i,\wbb)=\th^n(D^@,\qbb_j,\wbb)\Big)$$
where $\vvq_i:=\bbq_i\cap D$ and $\qbb_i:=\bbq_i\cap D^@$.

For $i<j$ in $T^1$ denote 
$$[\bbq_i,\bbq_j]^J_a:=[\vvq_i,\vvq_j]^J_a\cup\qbb_i$$
and note that by the preservation theorem the results of the shufflings 
are representatives for elements of\/ $\G$ i.e. 
$$i,j\in T^1 \imp C\mo U([\bbq_i,\bbq_j]_a^J,\bbw).$$
Define $k^*(i,j)$ for $i<j$ from  $T^1$ by
$$k^*(i,j) := \min\Big\{ k : \big(k\in T^1\ \&\
[\bbq_i,\bbq_j]_a^J\sim\bbq_k\big)\vee\big(k=N^*_5\big)\Big\}.$$
There is a subset $T\sb T^1$ of size $N^*_6$, large enough, \st for all 
$\bbq_i,\bbq_j,\bbq_l$ with $i<j<l$ from $T$ the following five statements 
have a constant truth value: $k^*(j,l)=i$, $k^*(i,l)=j$, $k^*(i,j)=i$, 
$k^*(i,j)=j$, $k^*(i,j)=l$.
Moreover, as usual if there are $i<j$ in $T$ with $k^*(i,j)\in T$ then
either for every $i<j$ in $T$, $k^*(i,j)=i$ or for every $i<j$
in $T$, $k^*(i,j)=j$.

If there isn't such a pair choose $\bbv_T$ \st 
$$\bigwedge_{i\in T}[C\mo R(\bbq_i,\bbv_T,\bbw)]\wedge
\bigwedge_{i<j\in T}[C\mo \neg R([\bbq_i,\bbq_j]^J_a,\bbv_T,\bbw)]$$
as $N^*_6$ is big enough there are $i<j$ from $T$ with:

\centerline{$\ath^{m(*)}(D;\vvq_i,\bbv_T\cap D,\vvw)=
\ath^{m(*)}(D;\vvq_j,\bbv_T\cap D,\vvw)$}

\centerline{$\th^n(D^@;\qbb_i,\bbv_T\cap D^@,\wbb)=
\th^n(D^@;\qbb_j,\bbv_T\cap D^@,\wbb)$}

\\By the preservation theorem and the composition theorem we get 
$$ \ \ (*) \ \
\th^n(C;[\bbq_i,\bbq_j]_a^J,\bbv_T,\bbw)=\th^n(C;\bbq_i,\bbv_T,\bbw).$$
Therefore $C\mo R([\bbq_i,\bbq_j]_a^J,\bbv_T,\bbw)$ which is a contradition. 

It follows: either $i,j\in T \imp [\bbq_i,\bbq_j]_a^J\sim\bbq_i$ or 
$i<j\in T \imp [\bbq_i,\bbq_j]_a^J\sim\bbq_j$.

Now choose $i,j\in T$ \st 

\centerline{$\ath^{m(*)}(D;\vvp_i,\vvp_{{2N^*_4}-i},\vvq_i,\vvw) =
\ath^{m(*)}(D;\vvp_j,\vvp_{{2N^*_4}-j},\vvq_j,\vvw)$}

\centerline{$\th^n(D^@;\pbb_i,\pbb_{{2N^*_4}-i},\qbb_i,\wbb) =
\th^n(D^@;\pbb_j,\pbb_{{2N^*_4}-j},\qbb_j,\wbb)$}

\\Shuffle along $a$ and $J$ and get a contradiction as before to the 
definition of the $\bbq_i$'s.

\

\

We have proved the following:

\pt Theorem 7.3. Let 
$\lan U(\bbx,\bbz),E(\bbx,\bby,\bbz),R(\bbx,\bby,\bbz)\ran$ be a sequence of 
formulas of dimension $d$ and depth $n$. 

\\Then there is $K<\om$, that depends only on $d$ and $n$ such that, in $V^P$, 
for no chain $C$ and parameters $\bbw\sb C$:

$(i)$ \/ $\II=\lan U(\bbx,\bbw),E(\bbx,\bby,\bbw),R(\bbx,\bby,\bbw)\ran$
is an interpretation for some $\G\in\Ga_k$ in $C$,

$(ii)$ \/ $D$, the minimal $(K_1,K_2)$-major initial (or final) segment 
for $\II$, is isomorphic to a regular cardinal $\la>\ale$.
\vv\qed

\

\bf Reduced Shufflings: \rm \ There are two main difficulties that face us in 
the general context. The first one is that the preservation theorem is 
formulated only in the context of well ordered chains. We can try and solve 
this by choosing a cofinal sequence through the chain and shuffle along this 
sequence. However the second difficulty is that a \sc\ that has the 
cardinality of\/ $\cf(D)$ (where $D$ is the minimal major initial segment)   
can't be generic with respect to subsets of\/ $D$ when $|D|>\cf(D)$. 
The solution for both this difficulties lies in the observation that what we 
really shuffle are not subsets of the chain but rather partial theories.

Suppose that we are given a chain $C$, with $\cf(C)=\la>\ale$ and some  
$\bb A\sb C$ of length $l$. For simplicity we assume that the chains have a 
first element $\min(C)$. Choosing a cofinal sequence $E=\lan\al_i:i<\la\ran$ 
in $C$ \st $\al_0=\min(C)$ and defining 
$s_i:=\th^n(C,\bb A)\res_{[\al_i,\al_{i+1})}$ we get by the composition 
theorem that 
$$\th^n(C,\bb A)=\sum_{i<\la}s_i.$$
Concentrating on the chain $(\la,<)$ we define a sequence 
$\bbh=\bbh_{\bba}=\lan P_t:t\in T_{n,l}\ran$ where for $t\in T_{n,l}$, 
$P_t:=\{i<\la:s_i=t\}$. By the Feferman-Vaught theorem (1.9) we know that 
$\th^n(C;\bba)$ is determined by $\th^m(\la;\bbh)$ where $m=m(n,l)$ depends 
only on $n$ and $l$. 

\

\pt Lemma 7.4. Let $C$ be a chain with cofinality $\la>\ale$ and let 
$n,l\in\NN$. Then, there are $m(*),l(*),\be(*)\in\NN$, all depending only on 
$n$ and $l$, \st

(a) there is a 1-1 function $\bbx\mapsto\bbh_{\bbx}$ \st for every $\bba\sb C$
of length $l$ there is $\bbh_{\bba}\sb\la$ of length $l(*)$ and 
$\th^n(C;\bba)$ is determined by 
$\th^{m(*)}(\la;\bbh_{\bba})$,

(b) $\be(*)$ codes a Turing machine that computes $\th^n(C;\bba)$ from 
$\th^{m(*)}(\la;\bbh_{\bba})$. 
\vv
\pd Proof. Choose a cofinal $E=\lan\al_i:i<\la\ran\sb\la$ ($\al_0=\min(C)$).  
Let $\bbh_{\bba}$ be as above $\l(*)=|T_{n,l}|$ and (a) is clear from the 
previous discussion. The computability in clause (b) is clear from the fact 
that $T_{m(*),l(*)}$ and $T_{n,l}$ are both finite. 
\vv\qed

\

\pd Remark. Of course we don't really lose generality by assuming that 
$C$ has a minimal element. If\/ $C$ interprets $\G$ by $\II$ and does'nt have 
one then we can always construct $C^*=C\cup\{-\infty\}$ and interpret $\G$ on 
$C^*$ by some $\II^*$ having the same depth and dimension $d+1$ (add 
$-\infty$ as a parameter). So instead of taking $K=K(n,d)$ we use 
$K=K(n,d+1)$ for getting a contradiction. 

\

The discussions above justify the following definition: 

\pd \dd 7.5. Let $n,d\in\NN$, and $\lan t_k:k<|T_{n,d}|\ran$ be the list of 
the possibilities $T_{n,d}$.

(1) $\T=(\la,\bbh)$ is a {\sl pre-chain} if\/ $\la>\ale$ is a regular cardinal, 
and $\bbh=\lan P_k:k<|T_{n,d}|\ran$ is a partition of\/ $\la$.

(2) We identify $(\la,\bbh)$ with $\lan s_i:i<\la\ran$ when  
$i\in P_k\iff s_i=t_k$.

(3) $\big((\la,\bbh),E\big)$ is a {\sl guess} for $(C,\bba)$ if:

\ \ \ $(i)$ $E=\lan\al_i:i<\la\ran\sb\la$ is cofinal in $C$ and $\al_0=\min(C)$,

\ \ \ $(ii)$ $\bba\sb C$ and $\lg(\bba)=d$,

\ \ \ $(iii)$ $\th^n(C;\bba)\res_{[{\al_i},{\al_{i+1}})}=s_i$ when
$(\la,\bbh)=\lan s_i:i<\la\ran$.
\vvv

\

Next we claim that the guesses (which are well ordered chains of the 
correct cardinality) represent faithfully the guessed chain.

\pd \dd 7.6. Suppose that

a. $C$ is a chain with cofinality $\la>\ale$,

b. $\bba,\bbb\sb C$ have length $d$,

c. $E=\lan\al_i:i<\la\ran$ is cofinal in $C$ and $\al_0=\min(C)$,

d. $J=\lan\be_j:j<\la\ran\sb\la$ is a club and $a\sb\la$ a \sc.

The {\sl reduced shuffling} of\/ $\bba$ and $\bbb$ along $E$, $J$ and $a$, 
denoted by $[\bba,\bbb]^{J,E}_a$ is defined by:

$$[\bba,\bbb]^{J,E}_a:=\bigcup_{j\in a}
\big(\bba\cap[\al_{\be_j},\al_{\be_{j+1}})\big)\cup
\bigcup_{j\not\in a}\big(\bbb\cap[\al_{\be_j},\al_{\be_{j+1}})\big)$$
\vvv

\

\pd Fact 7.7. If\/ $C$, $\bba$, $\bbb$, $J$ and $a$ are as above, 
$\big((\la,\bbh_{\bba}),E\big)$ a guess for $(C,\bba)$ and 
$\big((\la,\bbh_{\bbb}),E\big)$ a guess for $(C,\bbb)$ then 
$$[\bbh_{\bba},\bbh_{\bbb}]^J_a=\bbh_{[{\bba},{\bbb}]^{J,E}_a}$$
\vv
\pd Proof. Straightforward. 
\vv\qed

\

\pd \dd 7.8. For $C$, $\bba\sb C$, $E\sb C$ as above and $a\sb\la$ a \sc, 
define 
$$\ath_E^n(C;\bba):=\ath^n(\la,\bbh_{\bba})$$
\vvv

\

\pt Lemma 7.9. For every $n,d\in\NN$ there is $k(*)=k(n,d)\in\NN$ \st if 

1. $C$ is a chain and $\cf(C)=\la>\ale$,

2. $\bba,\bbb\sb C$ are of length $d$,

3. $E$ is cofinal in $C$,

4. $a\sb\la$ is a \sc,

\\then
$$\ath_E^{k(*)}(C;\bba)=\ath_E^{k(*)}(C;\bbb)\ \imp \
\th^n(C;\bba)=\th^n(C;\bbb)=\th^n(C;[\bba,\bbb]^{J,E}_a).$$
\vv
\pd Proof. Let $k(*)$ be $k\big(m(*),|T{n,d}\big|)$ where $m(*)$ is 
$m(*)(n,d)$ from the preservation theorem, and $k(\al,\be)$ is the 
``Feferman-Vaught'' number as in theorem 1.9. 
\vv\qed

\

\pt Lemma 7.10. Let $\lan U(\bbx,\bbz),E(\bbx,\bby,\bbz),R(\bbx,\bby,\bbz)\ran$ 
be a sequence of formulas of dimension $d$ and depth $n$. 

Then there is $K<\om$, that depends only on $d$ and $n$ \st in $V^P$,  
for no chain $C$ and parameters $\bbw\sb C$:

$(i)$ \ $\II=\lan U(\bbx,\bbw),E(\bbx,\bby,\bbw),R(\bbx,\bby,\bbw)\ran$
is an interpretation for some $\G\in\Ga_k$ in $C$,

$(ii)$ \ $D$, the minimal $(K_1,K_2)$-major initial (or final) segment 
for $\II$, satisfies $\cf(D)=\la>\ale$.
\vv

\pd Proof. We will follow the previous procedures, this time choosing $K$ 
big enough with respect to $k(*)$ as above and not $m(*)$ as usual. 
Assume that $\II$ is an interpretation of\/ $\G\in\Ga_K$ and suppose first that 
$C=D$. Fix $E\sb C$ cofinal, of order type $\la$. As $|E|=\la$, it belongs to 
the intermediate $V^{P_{\le\la}}$.

Let $\lan\bbu_i:i<|\G|\ran$ a list of the representatives and after the 
preliminary colouring we remain with a \sg list $B:=\lan\bbu_i:i\in S\ran$, 
($|S|=N_1$ big enough) having now the same $\ath^{k(*)}_E(C;\bbu_i)$.
Let $B_1=\lan\bbv_j:j<|S|^{N_2}\ran$ a list of the representatives 
for elements separating subsets of\/ $B$ of size $N_2$ from their complements.

Let $\lan\big((\la,\bbh_{\al,\be,\ga},E\big):\al,\be,\ga<N_1^{N_2}\ran$ 
be a list of all the guesses for chains of the form 

\\$(C;\bba_0,\bba_1,\bba_2,\bbw)$ with $\bba_i\in B\cup B_1$ for $i<3$.

Choose $J\sb\la$, a $k(*)$ suitable club for all the guesses, and a generic 
\sc\ $a\sb\la$. Start shuffling $\lan\bbu_i:i\in S\ran$ 
(i.e. the respective guesses).  A statement of the form $\bbu_\al\sim\bbu_\be$ 
is translated to ``$\th^{m(*)}(\la;\bbh_{\bbu_\al},\bbh_{\bbu_\be})$ is \st 
$C\mo E(\bbu_\al,\bbu_\be,\bbw)$''.

Repeating the usual steps we get $\lan\bbu_i:i\in A\ran$ \st w.l.o.g 
$[\bbu_i,\bbu_j]^{J,E}_a\sim\bbu_i$ for every $i<j$ from $A$.
Using genericity we can show also that $[\bbu_j,\bbu_i]^{J,E}_a\sim\bbu_i$
as well. 

Now choose a sequence of separating representatives $\lan\bbq_i:i<|A|/2\ran$ 
from $B_1$ above (so $J$ is suitable for them as well) and get a 
contradiction as usual. 

\

In the case $D\not=C$ we combine the above with the previous proof: 
the result of the shuffling of a pair of representatives $\bbu_\al$ and 
$\bbu_\be$ (coinciding outside $D$) is:

$\big\{$the result of the reduced shuffling of\/ $\bbu_\al\cap D$ and $\bbu_\be\cap D$ $\big\}$ 
$\bigcup$  $\big\{\bbu_\al\cap(C\sm D)\big\}$..

\\And we work in $D$. 
\vv\qed

\

As an $\om$-random graph is $K$-random for each $K<\om$ we proved:

\pt Theorem 7.11. In $V^P$:

A. If\/ $\lan \, C_K,\ I,\ \{\bbw_K : K\in A\}\;\ran$ is a uniform 
interpretation of\/ $\Ga_{\rm fin}$ in the monadic theory of order and 
$D_K\sb C_K$ are the minimal major initial (or final) segments of the 
interpretations, then $\cf(D_K)\le\ale$ for every large enough $K$.

B. If\/ $\II=\lan U(\bbx,\bbw),E(\bbx,\bby,\bbw),R(\bbx,\bby,\bbw)\ran$ is an 
interpretation for $RG_\om$ in a chain $C$ and $D\sb C$ is the 
minimal major initial (or final) segment then $\cf(D)\le\ale$.
\vv\qed

\

In the next section we will show that ``$\cf(D)\le\ale$'' is impossible 
even from ZFC.
\sec\sec
\hh 8. Short Chains                                                     \par 
Recall that a short chain is a chain that does not embed $(\om_1,<)$ and the 
inverse chain $(\om_1,>)$. Our aim in this section is to prove, from ZFC, 
the non exsistence of interpretations in short chains. In fact we show (and 
this is the only possibility when $C$ is short) the non exsistence of 
interpretations with $\cf(D)\le\ale$.  
 
\pd \dd 8.1. An interpretation $\II$ of\/ $\G\in\Ga_K$ in a chain $C$ is a 
{\sl short interpretation} if the minimal $(K_1,K_2)$-major initial (or final)
segment for $\II$, has cofinality $\ale$.
\vvv

\

The case $\cf(D)<\ale$ is impossible:

\pt Fact 8.2. Let $\II$ be an interpretation of some $\G\in\Ga_K$ in $C$. 
Let $D$ be the $(K_1,K_2)$-major initial segment. Then (if\/ $K$ is 
sufficiently big with respect to $d(\II)$ and $n(\II)$), $D$ does not have 
a last element.
\vv
\pd Proof. When $K$ is big enough we have $M(K,n,d)/{m(K,n,d)}>2$ (by 3.10) 
and this what we need. Now if\/ $D=D'\cup\{x\}$ where $x$ is the last element 
of\/ $D$ then, from the definitions, easily $\#(D)/{\#(D')}\le 2$. But 
$D'$ is minor and this is a contradiction.  
\vv\qed

\

\bf Assumptions. \rm From now on we are assuming towards a contradiction:

1. \  $\II=\lan U(\bbx,\bbw), E(\bbx,\bby,\bbw), R(\bbx,\bby,\bbw)$ is an 
interpretation for some $\G\in\Ga_K$ in a chain $C$. 
$n(\II)=n$ and $d(\II)=d$,

2. \ $K=K(n,d)$ is big enough (we will elaborate later),

3. \ $C$ has a minimal element (almost w.l.o.g by a previous remark),

4. \  $C$ is the minimal major initial segment for $\II$,

5. \ $\cf(C)=\ale$.

\

The next definition is the current replacement of\/ $m(*)$-suitable club:

\pd \dd 8.3. Let $\lan\bbu_i:i<i^*\ran$ be with $\bbu_i\sb C$, $\lg(\bbu_i)=d$. 
Let $E=\lan\al_k:k<\om\ran\sb C$ be increasing in $C$. 
$E$ is an {\sl $r$-suitable sequence for $\lan\bbu_i:i<i^*\ran$} if 

1. \ $E$ is cofinal in $C$ and $\al_0=\min(C)$,

2. \ For every $i<j<i^*$ there is $t_{i,j}\in T_{r,3d}$ \st for every 
$0<k<\om$: 

\centerline{$\th^r(C;\bbu_i,\bbu_j,\bbw)\res_{[{\al_0},{\al_k})}=t_{i,j}$,}

3. \ For every $i<j<i^*$ there is $s_{i,j}\in T_{r,3d}$ \st for every 
$0<k<l<\om$: 

\centerline{$\th^r(C;\bbu_i,\bbu_j,\bbw)\res_{[{\al_k},{\al_l})}=s_{i,j}$.}
\vvv

\

\\$r$-suitable sequences exist: 
\pt Claim 8.4. 

1. \ Suppose that $\bbu,\bbv\sb C$ are of length $d$ and 
$E=\lan\al_k:k<\om\ran$ is $r$-suitable for $\bbu,\bbv$. Let $E_1\sb E$ be 
infinite with $\al_0\in E_1$. Then $E_1$ is $r$-suitable for $\bbu,\bbv$.

2. \ Let $\bbu,\bbv\sb C$ be as above and let $E=\lan\al_k:k<\om\ran$ be 
cofinal with $\al_0=\min(C)$. Then there is $E_1\sb E$ that is
$r$-suitable for $\bbu,\bbv$.

3. \ For every finite family $\lan\bbu_i:i<i^*\ran$ 
with $\bbu_i\sb C$, $\lg(\bbu_i)=d$ there is an $r$-suitable $E\sb C$.
\vv
\pd Proof. The first part is immediate. For proving 2. let 
$\bbu$, $\bbv$, $E$ be given. Let $f\colon[\om\sm\{0\}]^2\to
|T_{r,3d}|\times|T_{r,3d}|$ be a colouring defined (for $0<k<l<\om$) by
$$f(k,l)=\big\lan\th^r(C;\bbu_i,\bbu_j,\bbw)\res_{[{\al_0},{\al_k})},
\th^r(C;\bbu_i,\bbu_j,\bbw)\res_{[{\al_k},{\al_l})}\big\ran$$
Let $u\sb\om$ be infinite, homogeneous with respect to $f$ ($T_{r,3d}$ 
is finite). Define $E_1:=\{\al_0\}\cup\{\al_k:k\in u\}$.

The third part is immediate by 1. and 2. 
\vv\qed

\

We will assume $\sqrt{K}\gg N_0\gg N_1\gg N_2\gg0$, all depending only on $n$ 
and $d$.

\

Let $M=|\G|$ and let $\lan\bbu_i:i<M\ran$ be a list of representatives for the 
elements of $\G$. Let $f\colon[M]^2\to|T_{n+d,3d}|$ be defined by 
$$f(i,j)=\th^{n+d}(C;\bbu_i,\bbu_j,\bbw).$$
We may assume that there is $S$ of size $N_0$, \sg with respect to 
$f$ and $(m+1)$, where $m$ is the bouquet size of minor segments.

Let $E=\lan\al_k:k<\om\ran\sb C$ be $(n+d)$-suitable for 
$\lan\bbu_i\ee\bbw:i\in S'\ran$, (by 8.4).
Let $S\sb S'$ be of size $N_1$ \st for every $i<j$ and $r<s$ from $S$, and 
for every $0<k<l<\om$:
$$\th^{n+d}(C;\bbu_i,\bbu_j,\bbw)\res_{[{\al_0},{\al_k})}=
\th^{n+d}(C;\bbu_r,\bbu_s,\bbw)\res_{[{\al_0},{\al_k})}:=t$$
and 
$$\th^{n+d}(C;\bbu_i,\bbu_j,\bbw)\res_{[{\al_k},{\al_l})}=
\th^{n+d}(C;\bbu_r,\bbu_s,\bbw)\res_{[{\al_k},{\al_l})}:=s.$$
This is possible by the definition of\/ $(n+d)$-suitability (and as $N_0$ 
is big enough). By the composition theorem for every $i<j$ in $S$:
$$\th^{n+d}(C;\bbu_i,\bbu_j,\bbw)=t+\sum_{k<\om}s\,.$$

\

\pd \dd 8.5. For $u\sb\om$ define the {\sl shuffling of\/ $\bbu_i$ and 
$\bbu_j$ along $u$} by 
$$[\bbu_i,\bbu_j]_u:=\bigcup_{k\in
u}\big(\bbu_i\cap[\al_k,\al_{k+1}\big)\cup\bigcup_{k\not\in u}
\big(\bbu_j\cap[\al_k,\al_{k+1}\big)$$
\vvv
\

\pt Claim 8.6. For every $i<j$ in $S$, for every $u\sb\om$, 
$C\mo U([\bbu_i,\bbu_j]_u,\bbw)$.
\vv
\pd Proof. By suitability of\/ $E$ and definition of\/ $S$ there are $t_0$ and
$s_0$ \st for every $i\in S$:

$(i)$ \  $\th^n(C;\bbu_i,\bbw)\res_{[{\al_0},{\al_k})}=t_0$ \ for every $0<k<\om$,

$(ii)$ \  $\th^n(C;\bbu_i,\bbw)\res_{[{\al_k},{\al_l})}=s_0$ \ for every $0<k<l<\om$,

$(iii)$ \ $\th^n(C;\bbu_i,\bbw)=t_0+\sum_{k<\om}s_0$.

\\By the definition of shuffling, for every $u\sb\om$ and $i<j$ in $S$,
$$\th^n(C;[\bbu_i,\bbu_j]_u,\bbw)=t_0+\sum_{k<\om}s_0=\th^n(C;\bbu_i,\bbw).$$
Therefore $C\mo U([\bbu_i,\bbu_j]_u,\bbw)$.
\vv\qed

\

Define now:

$e:=\{2k:k<\om\}$,

$o:=\{2k+1:k<\om\}$,

$p:=\om\sm\{0\}$,

$q:=\{0\}$.

\

Let 
$$k(i,j) := \min\Big\{k:\big(k\in S\ \&\ 
[\bbu_i,\bbu_j]_e\sim\bbu_k\big)\vee\big(k=|\G|\big)\Big\}.$$
By bigness of\/ $N_1$ there is $A\sb S$ of size $N_2$ \st for every 
$\bb U_i$, $\bb U_j$, $\bb U_l$ with $i<j<l$ from $A$, the following, usual, 
five statements have the same truth value:

$k(j,l)=i$,

$k(i,l)=j$,

$k(i,j)=i$,

$k(i,j)=j$,

$k(i,j)=l$.

\\Moreover, (the usual proof) if for some $i<j$ in $A$, $k(i,j)\in A$ then:
either for every $i<j$ in $A$, $k(i,j)=i$ or for every $i<j$ in $A$, 
$k(i,j)=j$.

\

Let's find $i<j$ in $A$ with $k(i,j)\in A$: if we can't then there is 
some $\bbv_A$ that separates between
$\A_1:=\{\bbu_i:i\in A\}$ and $\A_2:=\{\bbu_l:(\ex i<j\in
A)([\bbu_i,\bbu_j]_e\sim\bbu_l\}$. i.e. 
$$\bigwedge_{\bbu_i\in \A_1}\Big(C\mo R(\bbu_i,\bbv_A,\bbw)\Big)\wedge
\bigwedge_{{\bbu_i}\in{\A_2}}\Big(C\mo\neg R(\bbu_i,\bbv_A,\bbw)\Big).$$
We may assume that $E$ is suitable also for $\bbv_A$ (there are finitely 
many possibilities for $\bbv_A$ after choosing $\lan\bbu_i:i\in S'\ran$).  
As $N_2$ is big enough there are  $i<j$ in $A$ \st for every $0<k<l<\om$ 
$$\th^n(C;\bbu_i,\bbv_A,\bbw)\res_{[{\al_0},{\al_k})}=
\th^n(C;\bbu_j,\bbv_A,\bbw)\res_{[{\al_0},{\al_k})}$$
and
$$\th^n(C;\bbu_i,\bbv_A,\bbw)\res_{[{\al_k},{\al_l})}=
\th^n(C;\bbu_j,\bbv_A,\bbw)\res_{[{\al_k},{\al_l})}.$$
It follows that
$$\th^n(C;\bbu_i,\bbv_A,\bbw)=
\th^n(C;[\bbu_i,\bbu_j]_e,\bbv_A,\bbw)$$
and $\A_1\cap\A_2\ne\em$, a contradiction. We conclude, 
$$(*) \ \ (\ex i<j){\rm in}\ A\ {\rm such\ that}\ 
[\bbu_i,\bbu_j]_e\sim\bbu_i\ or\ [\bbu_i,\bbu_j]_e\sim\bbu_j.$$

\

\pt Fact 8.7. For every $i,j$ in $A$ (in fact in $S'$):
$[\bbu_i,\bbu_j]_q\sim\bbu_j$ and $[\bbu_i,\bbu_j]_p\sim\bbu_i$.
\vv
\pd Proof. Let's prove the first statement (the second is proved similarly).
By claim 8.6, $[\bbu_i,\bbu_j]_q$ is a representative hence is equivalent to 
$\bbu_l$ for some $l<|\G|$. Suppose that $l>j$. By semi-homogeneity of 
$S'$ (therefore of\/ $A$) there are $j<l_0<l_1\ldots<l_{m+1}$ \st 
$$\bigwedge_{r<m+1}\th^{n+d}(C;\bbu_j,\bbu_{l_r})=\th^{n+d}(C;\bbu_j,\bbu_l).$$
By definition of\/ $q$, $\bbu_l$ belongs to the vicinity of\/ $\bbu_j$. 
As ``belonging to the vicinity'' is determined by $\th^{n+d}$ we get $m+1$ 
pairwise nonequivalent representatives in $[\bbu_j]$. This is 
impossible by lemma 6.4. The same holds if we assume $l<j$. Therefore we must 
conclude $l=j$ i.e. $[\bbu_i,\bbu_j]_q\sim\bbu_j$.
\vv\qed

\

Returning to the representatives $\bbu_i$ and $\bbu_j$ we got in $(*)$ above,
suppose first that $[\bbu_i,\bbu_j]_e\sim\bbu_i$. We will show that

(1) \ $[\bbu_i,\bbu_j]_e\sim\bbu_i \imp [\bbu_i,\bbu_j]_o\sim
[\bbu_i,\bbu_j]_q$,

(2) \  $[\bbu_i,\bbu_j]_e\sim [\bbu_i,\bbu_j]_o$.

\\It will follow that $[\bbu_i,\bbu_j]_e\sim [\bbu_i,\bbu_j]_q$ and by the 
previous fact $\bbu_i\sim\bbu_j$ which is a contradiction.

\

For showing (1) it is enough to show that
$$\th^n(C;\bbu_i,[\bbu_i,\bbu_j]_e,\bbw)=
\th^n(C;[\bbu_i,\bbu_j]_o,[\bbu_i,\bbu_j]_q,\bbw).$$
Remembering how $S$ was chosen we get 

\

\\$\th^n(C;\bbu_i,[\bbu_i,\bbu_j]_e,\bbw)\res_{[\al_0,\al_1)}=
\th^n(C;\bbu_i,\bbu_i,\bbw)\res_{[\al_0,\al_1)}=$

\\$\th^n(C;[\bbu_i,\bbu_j]_o,[\bbu_i,\bbu_j]_q,\bbw)\res_{[\al_0,\al_1)}$.

\

\\$\th^n(C;\bbu_i,[\bbu_i,\bbu_j]_e,\bbw)\res_{[\al_{2k},\al_{2k+1})}=
\th^n(C;\bbu_i,\bbu_i,\bbw)\res_{[\al_{2k},\al_{2k+1})}=$

\\$\th^n(C;\bbu_j,\bbu_j,\bbw)\res_{[\al_{2k},\al_{2k+1})}=
\th^n(C;[\bbu_i,\bbu_j]_o,[\bbu_i,\bbu_j]_q,\bbw)\res_{[\al_{2k},\al_{2k+1})}$

\

\\$\th^n(C;\bbu_i,[\bbu_i,\bbu_j]_e,\bbw)\res_{[\al_{2k+1},\al_{2k+2})}=
\th^n(C;\bbu_i,\bbu_j,\bbw)\res_{[\al_{2k+1},\al_{2k+2})}=$

\\$\th^n(C;[\bbu_i,\bbu_j]_o,[\bbu_i,\bbu_j]_q,\bbw)\res_{[\al_{2k+1},\al_{2k+2})}$

\

and \ $\th^n(C;\bbu_i,[\bbu_i,\bbu_j]_e,\bbw)=
\th^n(C;[\bbu_i,\bbu_j]_o,[\bbu_i,\bbu_j]_q,\bbw)$ follows from the 
composition theorem.

\

For (2) note that

\\$\th^n(C;[\bbu_i,\bbu_j]_e,\bbu_i,\bbw)=$ \ 
$\th^n(C;[\bbu_i,\bbu_j]_e,\bbu_i,\bbw)\res_{[\al_0,\al_1)}+$

\\$\th^n(C;[\bbu_i,\bbu_j]_e,\bbu_i,\bbw)\res_{[\al_1,\al_2)}+
\th^n(C;[\bbu_i,\bbu_j]_e,\bbu_i,\bbw)\res_{[\al_2,\al_3)}+\ldots$

\\$=\th^n(C;\bbu_i,\bbu_i,\bbw)\res_{[\al_0,\al_1)}+
\th^n(C;\bbu_j,\bbu_i,\bbw)\res_{[\al_1,\al_2)}+
\th^n(C;\bbu_i,\bbu_i,\bbw)\res_{[\al_2,\al_3)}+$

\\$\th^n(C;\bbu_j,\bbu_i,\bbw)\res_{[\al_3,\al_4)}+
\th^n(C;\bbu_i,\bbu_i,\bbw)\res_{[\al_4,\al_5)}+\ldots$

\\and that

\\$\th^n(C;[\bbu_i,\bbu_j]_o,\bbu_i,\bbw)=$ \ 
$\th^n(C;[\bbu_i,\bbu_j]_o,\bbu_i,\bbw)\res_{[\al_0,\al_2)}+$ 

\\$\th^n(C;[\bbu_i,\bbu_j]_o,\bbu_i,\bbw)\res_{[\al_2,\al_3)}+
\th^n(C;[\bbu_i,\bbu_j]_o,\bbu_i,\bbw)\res_{[\al_3,\al_4)}+\ldots$

$=\th^n(C;\bbu_i,\bbu_i,\bbw)\res_{[\al_0,\al_2)}+
\th^n(C;\bbu_j,\bbu_i,\bbw)\res_{[\al_2,\al_3)}+
\th^n(C;\bbu_i,\bbu_i,\bbw)\res_{[\al_3,\al_4)}+$

$\th^n(C;\bbu_j,\bbu_i,\bbw)\res_{[\al_4,\al_5)}+
\th^n(C;\bbu_i,\bbu_i,\bbw)\res_{[\al_5,\al_6)}+\ldots$

\\By the composition theorem: 
$$\th^n(C;[\bbu_i,\bbu_j]_e,\bbu_i,\bbw)=\th^n(C;[\bbu_i,\bbu_j]_o,\bbu_i,\bbw).$$
That is:
$$[\bbu_i,\bbu_j]_e\sim\bbu_i\sim[\bbu_i,\bbu_j]_o.$$

\\Collecting the results we get:

$[\bbu_i,\bbu_j]_e\sim\bbu_i$ \ (this is the assumption),

$\bbu_i\sim[\bbu_i,\bbu_j]_o$ \ (by (2) above),

$[\bbu_i,\bbu_j]_o\sim[\bbu_i,\bbu_j]_q$ \ (by (1) above),

$[\bbu_i,\bbu_j]_q\sim\bbu_j$ \ (by fact 8.7).

\\Therefore, $\bbu_i\sim\bbu_j$ a contradiction.

\

\\We are therefore forced to assume that $[\bbu_i,\bbu_j]_e\sim\bbu_j$ but 
then we get the same way 
$\bbu_i\sim[\bbu_i,\bbu_j]_o$ (like (2) above), 
$[\bbu_i,\bbu_j]_o\sim[\bbu_i,\bbu_j]_p$ (like (1) above), 
$[\bbu_i,\bbu_j]_p\sim\bbu_j$ (by 8.7), and again $\bbu_i\sim\bbu_j$.

\

\

We assumed that $C$ is equal to $D$, the minimal major initial segment 
for simplicity. However, if  $D\ne C$ then following previous procedures we 
can easily chop off\/ $C\sm D$ and basically work inside $D$, getting a 
contradiction. 

So we have eliminated the possibilities that were left by theorem 7.11 and
proved:

\pt Theorem 8.8 (Non-Interpretability Theorem). There is a forcing notion $P$ 
such that in $V^P$ the following hold:  

\\(1) $RG_\om$ is not interpretable in the monadic theory of order.

\\(2) For every sequence of formulas $\II=\lan U(\bbx,\bbz), 
\ E(\bbx,\bby,\bbz), \ R(\bbx,\bby,\bbz)\ran$ there is $K^*<\om$, 
(effectively computable from $\II$), \st for no chain $C$, $\bb W\sb C$, and 
$K\ge K^*$ does $\lan U(\bbx,\bbw),$ $E(\bbx,\bby,\bbw),R(\bbx,\bby,\bbw)\ran$
interpret $RG_K$ in $C$. 

\\(3) The above propositions are provable in ZFC. if we restrict ourselves 
to the class of short chains. 
\vv\qed
\sec\sec 
\vfill\eject

\font\ba=cmr8
\font\bs=cmbxti10
\font\bib=cmtt12
\centerline{\bib REFERENCES}  

\

\\ \ {\bf [Ba]} \ J. B{\ba ALDWIN}, \ 
{\sl Definable second--order quantifiers}, \ 
{\bs Model Theoretic Logics}, \ 
(J. Barwise and S. Feferman, editors), 
\ Springer--Verlag, Berlin 1985, pp. 445--477.  \vvv
\\ \ {\bf [BaSh]} \ J. B{\ba ALDWIN} and S. S{\ba HELAH}, \ 
{\sl Classification of theories by second order quantifiers}, \ 
{\bs Notre Dame Journal of Formal Logic}, 
\ vol. 26 (1985) pp. 229--303.  \vvv
\\ \ {\bf [Gu]} \ Y. G{\ba UREVICH},  \ 
{\sl Monadic second--order theories}, \ 
{\bs Model Theoretic Logics}, \ 
(J. Barwise and S. Feferman, editors), 
\ Springer--Verlag, Berlin 1985, pp. 479--506.  \vvv
\\ \ {\bf [GMS]} \ Y. G{\ba UREVICH}, M. M{\ba AGIDOR} and S. S{\ba HELAH}, \ 
{\sl The monadic theory of\/ $\om_2$}, \ 
{\bs The Journal of Symbolic Logic}, 
\ vol. 48 (1983) pp. 387--398.   \vvv
\\ \ {\bf [GuSh]} \ Y. G{\ba UREVICH} and S. S{\ba HELAH}, \ 
{\sl On the strength of the interpretation method}, \ 
{\bs The Journal of Symbolic Logic}, 
\ vol. 54 (1989) pp. 305--323.   \vvv
\\ \ {\bf [GuSh1]} \ Y. G{\ba UREVICH} and S. S{\ba HELAH}, \ 
{\sl  Monadic theory of order and topology in ZFC}, \ 
{\bs Ann. Math. Logic}, 
\ vol. 23 (1982) pp. 179--182.   \vvv
\\ \ {\bf [GuSh2]} \ Y. G{\ba UREVICH} and S. S{\ba HELAH}, \ 
{\sl Interpreting second--order logic in the monadic theory of order}, \ 
{\bs The Journal of Symbolic Logic}, 
\ vol. 48 (1983) pp. 816--828.   \vvv
\\ \ {\bf [GuSh3]} \ Y. G{\ba UREVICH} and S. S{\ba HELAH}, \ 
{\sl The monadic theory and the `next world'}, \ 
{\bs Israel Journal of Mathematics}, 
\ vol. 49 (1984) pp. 55--68.   \vvv
\\ \ {\bf [LiSh]} \ S. L{\ba IFSCHES} and S. S{\ba HELAH}, \     
{\sl  Peano arithmetic may not be interpretable in the monadic theory of linear 
orders}, \ 
{\bs The Journal of Symbolic Logic}, 
\ accepted for publication.   \vvv
\\ \ {\bf [Sh]} \ S. S{\ba HELAH}, \ 
{\sl The monadic theory of order}, \ 
{\bs Annals of Mathematics}, 
\ ser. 2, vol. 102 (1975) pp. 379--419.   \vvv
\\ \ {\bf [Sh1]} \ S. S{\ba HELAH}, \ 
{\sl Notes on monadic logic, part B: complexity of linear orders}, \ 
{\bs Israel Journal of Mathematics}, 
\ vol. 69 (1990) pp. 99--116.   \vvv
\\ \ {\bf [TMR]} \ A. T{\ba ARSKI}, A. M{\ba OSTOWSKI} and R. M. R{\ba OBINSON}, \ 
{\bs Undecidable theories}, 
\ North-Holland Publishing Company, 1953, $xi$+98 pp.   \vvv
\end